\documentclass{article}
\usepackage[utf8]{inputenc}
\usepackage{rotating}
\usepackage{csquotes}
\title{Learning differential equation models from\\ stochastic agent-based model simulations}
\date{\today}

\author{John T. Nardini, Ruth E. Baker, Matthew J. Simpson, Kevin B. Flores}

\usepackage{natbib}
\usepackage{multirow}
\newcommand{\xmark}{\textcolor{red}{\ding{53}}}
\usepackage{amsmath}
\usepackage{geometry}
\usepackage{units}
\usepackage{amsfonts}
\usepackage{amssymb}
\usepackage{pifont}
\usepackage{multirow}
\usepackage{pdflscape}
\usepackage[dvipsnames]{xcolor}
\usepackage{bm}
\usepackage{hyperref}
\usepackage{ulem}

\definecolor{green}{RGB}{0,128,0}
\definecolor{red}{RGB}{200,0,0}
\newcommand{\checkmarkJTN}{\textcolor{green}{\ding{51}}}

\newcommand{\new}[1]{#1}
\newcommand{\old}[1]{}

\usepackage[ruled,vlined]{algorithm2e}

\usepackage{todonotes}

\newcommand{\mean}[1]{\langle #1 \rangle} 

\begin{document}

\maketitle

\begin{abstract}
    Agent-based models provide a flexible framework that is frequently used for modelling many biological systems, including cell migration, molecular dynamics, ecology, and epidemiology. Analysis of the model dynamics can be challenging due to their inherent stochasticity and heavy computational requirements. Common approaches to the analysis of agent-based models include extensive Monte Carlo simulation of the model or the derivation of coarse-grained differential equation models to predict the expected or averaged output from the agent-based model. Both of these approaches have limitations, however, as extensive computation of complex agent-based models may be infeasible, and coarse-grained differential equation models can fail to accurately describe model dynamics in certain parameter regimes. We propose that methods from the equation learning field provide a promising, novel, and unifying approach for agent-based model analysis. Equation learning is a recent field of research from data science that aims to infer differential equation models directly from data. We use this tutorial to review how methods from equation learning can be used to learn differential equation models from agent-based model simulations. We demonstrate that this framework is easy to use, requires few model simulations, and accurately predicts model dynamics in parameter regions where coarse-grained differential equation models fail to do so.  We highlight these advantages through several case studies involving two agent-based models that are broadly applicable to biological phenomena: a birth-death-migration model commonly used to explore cell biology experiments and a susceptible-infected-recovered model of infectious disease spread.
\end{abstract}

\section{Introduction}\label{sec:introduction}

Complex interactions between individuals are a crucial aspect of many biological processes: honeybees dance to direct others to food sources \citep{couzin_effective_2005}, cells push their neighbors to promote invasion during tumourigenesis \citep{schmidt_interstitial_2009}, and animal herds aggregate together to deter predation \citep{binny_living_2020}. Agent-based models (ABMs) are invaluable tools to simulate how such interactions between individuals scale to population-wide phenomena \citep{dorsogna_self-propelled_2006}. In an ABM, the states and decisions of individual agents are simulated using pre-defined rules to govern the agents' interactions and behaviour \citep{an_optimization_2017}. The ease of construction of ABMs by domain experts and modellers allows for complex models that can capture rich dynamical behaviour \citep{an_optimization_2017,mirams_chaste_2013}. 

\begin{figure}
    \centering
    \includegraphics[width=0.9\textwidth]{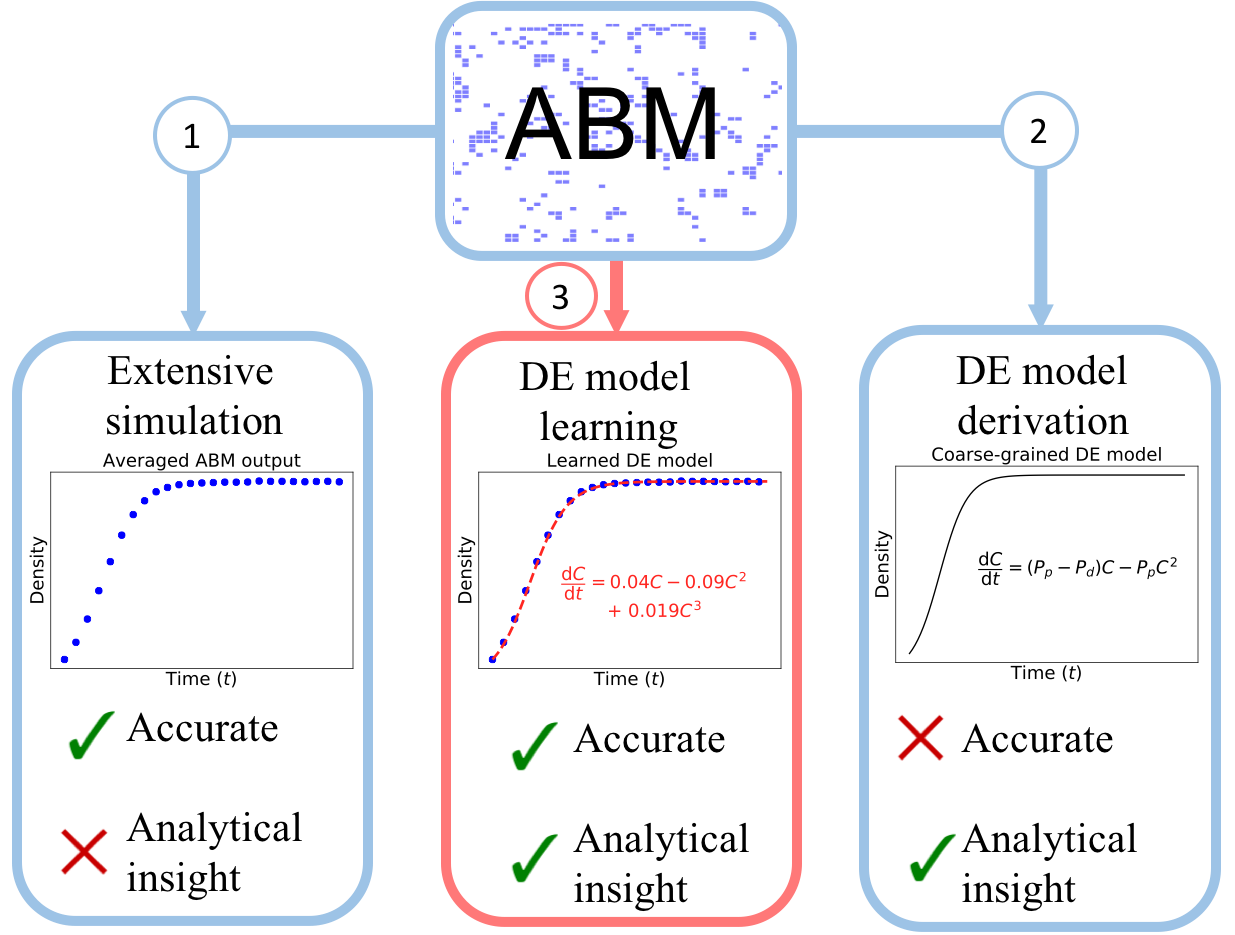}
    \caption{An illustration of current (blue) and proposed (red) methods to predict emergent ABM behaviour. Extensive simulation (Arrow 1) is performed by running many ABM simulations over a range of parameter values, and then using Monte Carlo techniques to average ABM output. While this approach will accurately predict ABM dynamics, it can be computationally intensive to perform. DE models derived using model coarse graining approaches (Arrow 2) can be analysed (e.g., using bifurcation analysis or perturbation methods). This technique is advantageous because such analytical methods do not require any computation. Unfortunately, coarse-grained models will provide inaccurate predictions in many parameter regimes. We propose that DE models can be learned from ABM simulation data using  techniques from equation learning (Arrow 3). This method is advantageous because it may only require a small number of ABM simulations, will lead to a DE model that can predict ABM dynamics accurately, and can be informed with analytical techniques.
}
    \label{fig:abm_eql_concept_figure}
\end{figure}

\old{Predicting the emergent behaviour of ABMs can be challenging, however, due to their inherent stochastic nature and analytic intractability \citep{irons_effect_2018}.}\new{There are many approaches to predicting the emergent behaviour of stochastic ABMs, each of which presents its own advantages and challenges.} \old{As such,} The most straightforward and commonly-used approach to interrogate ABMs is extensive Monte Carlo simulation using well-established computational algorithms \citep{an_optimization_2017} (Arrow 1 of Figure \ref{fig:abm_eql_concept_figure}). Average ABM behaviours at fixed parameter values can be inferred from many simulations from the central limit theorem \citep{seber_nonlinear_1988}, and the inverse problem of inferring model parameter distributions can be done with Markov chain Monte Carlo samplers \citep{beheshti_improving_2013}. Unfortunately, such extensive simulation of many ABMs may not be feasible due to significant computational costs involved. An alternative method to predict the emergent behaviour of an ABM consists of deriving \old{coarse-grained} differential equation (DE) models to approximate ABM output (Arrow 2 of Figure \ref{fig:abm_eql_concept_figure}). \new{Each ABM has a master equation that can be derived directly from the model rules. There are many approaches to simplify this master equation and approximate its dynamics with more tractable DE models.} The most commonly-used DE model approximations for ABMs are \emph{mean-field} models \citep{baker_correcting_2010,chaplain_bridging_2020,cruz_stochastic_2016,fadai_accurate_2019}. \new{Alternative formulations to mean-field models are also possible, see \citep{schnoerr_approximation_2017} for an extensive tutorial}. Mean-field models describe the evolution of population density over time (and possibly space) and can be derived by approximating  agent-agent interactions with locally-averaged agent densities \citep{fadai_accurate_2019}. \new{Mean-field} DE models are often simple to solve (either analytically or numerically), so they provide an advantageous alternative to extensive simulation of the ABM. Furthermore, such DE models are amenable to analytical techniques (including bifurcation, travelling wave, perturbation analysis), which can be used to predict how ABM output will change in response to variations in parameter values \citep{j_d_murray_asymptotic_2012,murray_mathematical_2002}.

Previous ecological studies have demonstrated some of the advantages of both extensive simulation and model coarse graining for ABM analysis \citep{bernoff_nonlocal_2013}. For example, Bernoff et al. \citep{bernoff_agent-based_2020} model the foraging behaviour of the Australian plaque locust with a discrete and stochastic ABM. In the model, individual locusts forage and feed on a given resource (representative of food) and, in turn, create a spatial gradient of this resource. The model robustly shows that individual locust behaviour drives the formation of this resource gradient and, in turn, determines how the averaged profile of locust density migrates and forms over time. The authors derive the mean-field partial differential equation (PDE) model  for this ABM and perform a travelling wave analysis to quantify how the locust population's invasion speed depends on the total mass of locusts. The mean-field PDE model is shown to match the ABM output well in biologically consistent parameter regimes. In addition, non-mean-field models have been considered to approximate other ABMs of locust behaviour. For example, the ABMs in \citep{dkhili_self-organized_2017} describe self-organizing  locust behaviours through rules governing locust attraction, repulsion, and alignment during foraging and invasion. By simulating the ABM over many different parameter values, Dkhili et al. \citep{dkhili_self-organized_2017} discovered three distinct population patterns (spot, band, and ribbon formations). Topaz et al. \citep{topaz_locust_2012} analysed a continuous  partial integro-differential equation as a representation of this locust flocking behaviour and used a linear stability analysis to provide analytical insights into which parameter values governing agent interactions lead to the formation of such spatial patterns. There are thus many scenarios in which DE models supplement ABM simulations to aid in our understanding of emergent behaviour.

Despite their wide use, coarse-grained models can provide misleading predictions of ABM dynamics in regions of parameter space in which the assumptions made during the coarse graining process do not hold \citep{baker_correcting_2010,fadai_accurate_2019,matsiaka_continuum_2017}.
Furthermore, it can be challenging to determine informative DE models for more complex ABMs. As one such example, Gallaher et al. \citep{gallaher_evolution_2013} constructed an ABM in which thousands of cells with different phenotypes compete for space during tumour growth. Each agent in the simulation is given an internal set of dynamic traits dictating how fast the agent moves in space and how frequently it divides. The intricate dynamics of this model allow for interesting findings of biological relevance, including how transmission of proliferation rates from parent to daughter cells alters the final trait landscape of the population and, in turn, the eventual physical clustering of the population. Formulating a DE model for this process from ABM rules would be challenging, however, due to the many different cell phenotypes and complicated rules between such cells.  Instead, methods to directly infer DE models from a small number of ABM simulations may provide a useful tool for modellers to determine the salient features necessary for modelling complex ABM dynamics. ABMs with evolving trait landscapes are becoming increasingly common to study tumour dynamics, so such learned DE models will be widely applicable to this growing field of research \citep{west_evolutionary_2016}.

\emph{Equation learning} (EQL) is a recent field from data science that aims to infer the dynamical systems model that best describes a given dataset \citep{brunton2016discovering}. The learned models can, in turn, be used to understand the system under study by providing a mechanistic description of observed dynamics or predicting how dynamics will change in response to different conditions. There has been much progress in this field over the past five years largely thanks to increases in computational power, and many EQL methods can accurately recover DE models from artificially simulated noisy data from DE models \citep{lagergren_learning_2020,rudy_data-driven_2017,zhang2018robust,zhang_robust_2019}. There have been some recent studies showing that equations can be learned from noisy experimental data \citep{lagergren_biologically-informed_2020}, and the EQL field is now in a position where EQL can, in principle, be used to aid in the development of DE models to approximate the dynamics of complex ABMs (Arrow 3 of Figure \ref{fig:abm_eql_concept_figure}). \new{In this review article, we will detail how a commonly-used EQL methodology \citep{brunton2016discovering} can be used to learn DE models that accurately describe ABM dynamics. In doing so, we also explore how EQL provides insight into ABM behaviours when traditional modelling approaches (\emph{e.g.}, coarse-graining) fail to capture ABM complexity, as well as EQL performance in practical situations, such as those with sparse data samples.}

Similar to coarse-grained DE models, learned DE models can be analysed using both computational algorithms and analytical techniques to infer the emergent behaviour that results from a given set of ABM, or estimate mechanistic parameters from data. Furthermore, such learned equations may have fewer computational requirements than ABMs if they can learn ABM dynamics from a small number of simulations. As a result, EQL methods may be tractable for learning DE models for a broad range of ABM simulations. There are many ways in which learned DE models may be useful for ABM analysis, including the discovery of novel DE models, predicting unobserved dynamics from complex ABMs, and enabling accurate parameter estimation.  Furthermore, as ABMs imitate many key features of biological systems (including stochasticity and heterogeneity), inferring DE models from ABM data is an intermediate step towards developing algorithms to aid the discovery of models from experimental data. 


This article is intended to serve as a review and tutorial on three separate but synergistic methods (extensive simulation, coarse-grained model derivation, and EQL) to infer the emergent behaviour of ABMs. The first two are frequently used for ABM analysis, and we propose that \new{analysis of a learned DE model from ABM data provides increased understanding of the mechanisms driving observed behaviour}.  We will showcase each of these methods and highlight the advantages and limitations of each. In Section 2, we discuss ABM setup and implementation as well as how simple ABM rules can be coarse grained directly into DE models. In Section 3, we discuss how methods from EQL can be used to learn DE models from ABM data directly. \new{Our goal in Section 4 is to present six questions (Q1 - Q6) relating to how EQL can aid in the analysis of ABM behaviour, and we address each question with case study examples.} We use two representative ABMs throughout: a birth-death-migration (BDM) model of population dynamics \citep{baker_correcting_2010}, and a susceptible-infected-recovered (SIR) model of infectious disease dynamics. We note, however, that the approaches discussed within this tutorial are broadly applicable to many social and biological phenomena that have been modelled by ABMs previously \citep{an_optimization_2017,bonabeau_agent-based_2002,cosgrove_agent-based_2019,interian_tumor_2017,stevens_stochastic_2000,stevens_aggregation_1997}. We make final conclusions, summarize the advantages and limitations of each method, and suggest future avenues for research in Section \ref{sec:conclusions}. Python code for all tutorials and case studies shown in this study are publicly available at  \href{https://github.com/johnnardini/Learning-DE-models-from-stochastic-ABMs}{{\color{blue} https://github.com/johnnardini/Learning-DE-models-from-stochastic-ABMs}}.

\section{coarse graining ABMS into DE models} \label{sec:ABM_deriv}

Coarse-grained DE models are now frequently used to investigate how rules governing individual behaviours translate to emergent behaviour at the population level \citep{anguige_one-dimensional_2009,baker_correcting_2010,fadai_accurate_2019}.  In this section, we illustrate this approach by introducing two simple ABMs and coarse graining them to give DE models. We consider two ABMs: (i) a BDM process in Section \ref{sec:BMD_process}, and (ii) a SIR model in Section \ref{sec:SIR_ABM}. While we focus on ordinary differential equation (ODE) models throughout this article, the derivation of PDE models in the presence of spatial heterogeneity can also be performed using extensions of the methods presented herein, as discussed in \citep{fadai_accurate_2019,johnston_how_2014,simpson_multi-species_2009}.

\subsection{A birth, death, and migration ABM}\label{sec:BMD_process}

We consider the BDM process, a lattice-based ABM in which agents are able to give birth, die, and move \citep{baker_correcting_2010}. This ABM is representative of many biological phenomena, for example, agents may represent cells during the wound healing process \citep{johnston_how_2014} or the invasion of animals in ecology \citep{bernoff_nonlocal_2013}. We begin by introducing the ABM rules in Section \ref{subsubsec:ABM_rules} and then coarse grain these rules into DE models and compare to ABM output in Section \ref{sec:ex1_logistic_deriv}.

\subsubsection{ABM Rules} \label{subsubsec:ABM_rules}
We use a two-dimensional square lattice with a lattice spacing of $\Delta$. We arbitrarily set $\Delta=1$ and assume the lattice has $X$ lattice sites in each spatial dimension. Each lattice site is indexed by $\alpha=(i,j)\in\mathbb{N}^2,\ i,j=1,\dots,X$. For each interior lattice site, $\alpha$, we define its neighbouring sites, $\mathcal{B}(\alpha)$, using the Von Neumann neighborhood $\mathcal{B} (\alpha)=\{(i,j+1),(i,j-1),(i+1,j),(i-1,j)\}$ and adjust this definition at boundary sites to enforce no-flux conditions. We designate the occupancy of each lattice site $\alpha$ as $\mu_\alpha(t) = 0$ if $\alpha$ is unoccupied at time $t$ or by $\mu_\alpha(t) = A$ if $\alpha$ is occupied by an agent at time $t$. For simplicity, we will also use the notation $0_\alpha(t)$ and $A_\alpha(t)$ when $\alpha$ is unoccupied and occupied, respectively.  To represent volume exclusion, or crowding, we assume that each lattice site can be occupied by a maximum of one agent at a time.


We next define how agents proliferate, migrate, and die. For proliferation events, we assume that agents proliferate with rate $P_p$ (formally, an agent will attempt to proliferate over an infinitesimal timestep of duration $\text{d} t$ with probability $P_p\text{d} t$). If an agent chooses to proliferate, then it will attempt to place its daughter cell into a neighbouring site $\beta\in\mathcal{B}(\alpha)$ with the choice of $\beta$ made uniformly at random. If the chosen site is occupied, the event is aborted. This process may be written as a bimolecular reaction with rate $P_p/4$: 
\begin{equation}
    A_\alpha + 0_\beta \xrightarrow{P_p/4} A_\alpha + A_\beta, \ \beta\in\mathcal{B}(\alpha). \label{eq:prolife_rate}
\end{equation}
The reaction rate here is divided by four because proliferating agents randomly choose one of their four neighbouring sites to place their daughter cell into. For death events, lattice sites  transition from occupied to unoccupied without any explicit crowding effects. We assume agents die with rate $P_d$, and write death events as a monomolecular reaction with rate $P_d$:
\begin{equation}
    A_\alpha \xrightarrow{P_d} 0_\alpha. \label{eq:death_rate}
\end{equation}
 Agents attempt to move with rate $P_m$. During migration events, agents randomly choose a neighbouring lattice site $\beta\in\mathcal{B}(\alpha)$ to attempt to move to. If the chosen site is already occupied, then the migration event is aborted. This process may be written as a bimolecular reaction with rate $P_m/4$:
\begin{equation}
    A_\alpha + 0_\beta \xrightarrow{P_m/4} 0_\alpha + A_\beta, \ \beta\in\mathcal{B}(\alpha).\label{eq:migrate_rate}
\end{equation}
The reaction rate here is divided by four because migrating agents randomly choose one of their four neighbouring sites to place their daughter cell into. \new{Following previous ABM studies, we have chosen to include empty space as an interacting agent in Equations (1)-(2) to incorporate the effects of volume exclusion. This choice converts migration into a bimolecular reaction instead of a monomolecular reaction \citep{baker_correcting_2010,cianci_molecular_2016,jin_extended_2018,johnston_mean-field_2012,mckane_stochastic_2004,simpson_multi-species_2009}. }

\begin{figure}
    \centering
    \includegraphics[width=.99\textwidth]{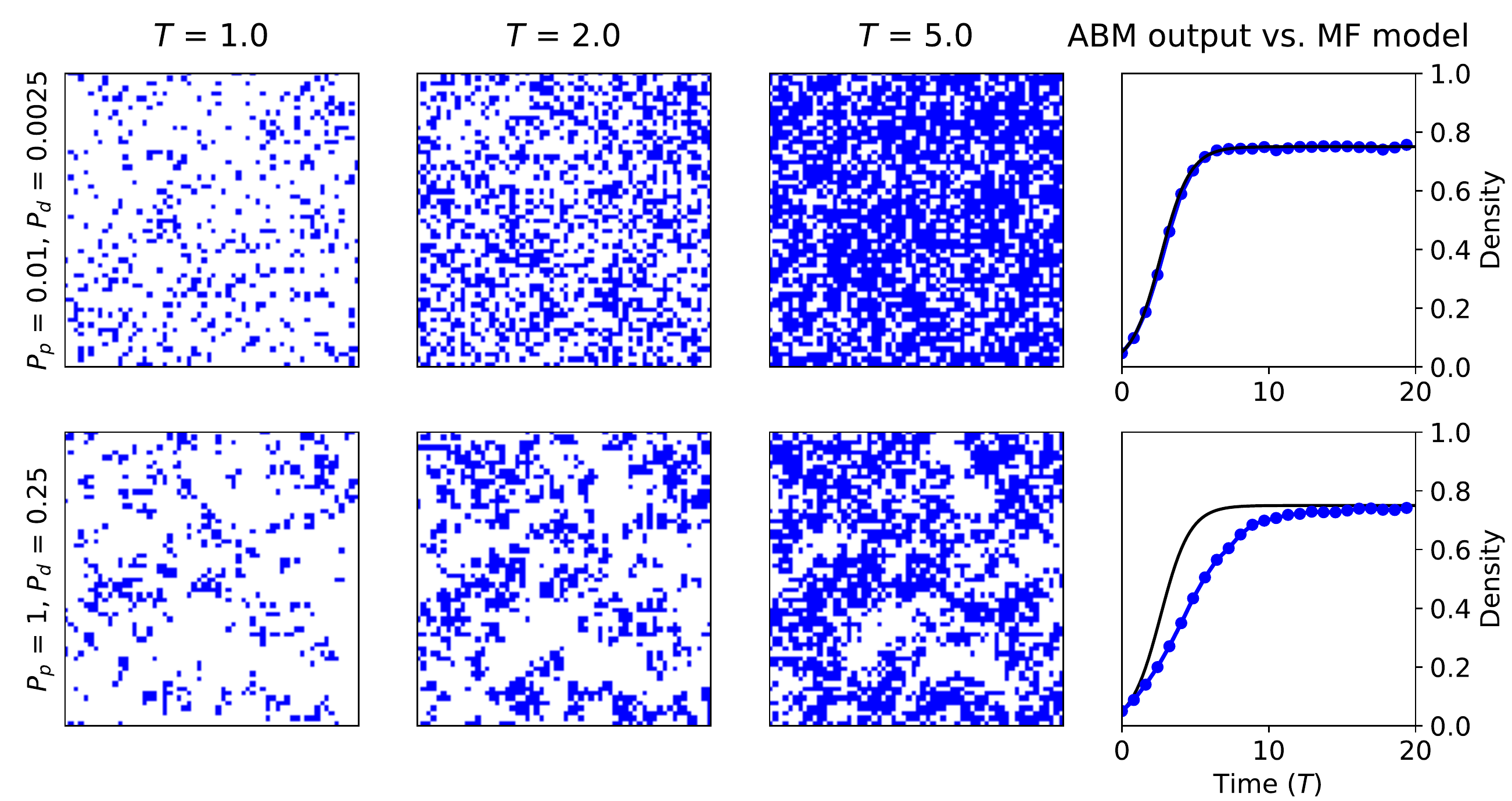}
    \caption{ABM simulation snap shots for the BDM ABM. Blue pixels denote $A_\alpha(t)$ and white pixels denote $0_\alpha(t)$. Simulations were computed with agent migration rate $P_m = 1$. The right-most column depicts \new{the output agent density from one ABM simulation}  (blue line and dots) against the solution of the logistic equation \eqref{eq:log_analytic_main} (solid black line). Quantities are plotted against nondimensionalised time $T=(P_p-P_d)t$ for ease of interpretation. The ABM was computed on a square lattice of length $X=120$.}
    \label{fig:ABM_sims}
\end{figure}

The BDM model can be simulated using the Gillespie algorithm \citep{gillespie_exact_1977}, which is \new{provided for the BDM process in Algorithm \ref{algo:gillespie} in Appendix \ref{sec:Gillespie}} \old{is summarised in \citep{fadai_accurate_2019}}. Each ABM simulation is initialized by placing agents uniformly at random throughout the lattice so that 5\% of lattice sites are occupied. \new{Note that each ABM simulation begins with a new initial configuration of agent locations throughout the lattice}. Reflecting boundary conditions are used at the boundaries of the lattice, which enforces a no-flux condition in the spatial domain. We use the following notation to summarize the output from an ABM simulation. To estimate the total agent density from the $n^\text{th}$ of $N$ identically prepared ABM simulations we compute
\begin{equation}
    C_{\text{ABM}}^{(n)}(t)= \dfrac{C^{(n)}(t)}{X^2}, \label{eq:ABM_dens}
\end{equation}
where $C^{(n)}(t)$ is the number of occupied sites at time $t$. To estimate the averaged agent density over time from $N$ identically prepared ABM simulations, we compute
\begin{equation}
    \mean{C_{\text{ABM}}(t)}= \dfrac{1}{N}\sum_{n=1}^N C_\text{ABM}^{(n)}(t). \label{eq:ABM_dens_avg}
\end{equation}
We depict snapshots of two simulations of the BDM process in Figure \ref{fig:ABM_sims}. The blue dots in the right-most column correspond to $C_\text{ABM}^{(1)}(t)$ from these individual simulations.

\subsubsection{Model coarse graining} \label{subsec:k-point}\label{sec:ex1_logistic_deriv}

ABM rules are often coarse grained into continuous DE models to aid in their analysis. Many previous studies \citep{baker_correcting_2010,fadai_accurate_2019,simpson_experimental_2013} have demonstrated that the mean-field DE model for the BDM process described in Section \ref{sec:BMD_process} is given by the logistic DE model 
\begin{equation}
    \dfrac{\text{d}}{\text{d}t}C(t) = P_pC(t)\big(1-C(t)\big) - P_dC(t). \label{eq:logistic}
\end{equation}
This model is advantageous in that it can be solved analytically to give
\begin{equation}
    C(t) = \dfrac{KC(0)\text{e}^{rt}}{K + C(0)(\text{e}^{rt}-1)},\label{eq:log_analytic_main}
\end{equation}
where $r = P_p - P_d$,  $K = (P_p-P_d)/P_p$, and $C(0)$ denotes the initial condition. The full derivation of this model is provided in Appendix A. To arrive at Equation \eqref{eq:logistic}, we use the \emph{mean-field assumption} that the occupancies of neighbouring lattice sites are independent, i.e. for all sites $\alpha$, $\mathbb{P}[A_\alpha(t),A_\beta(t)]=\mathbb{P}[A_\alpha(t)]\mathbb{P}[A_\beta(t)],\ \beta\in\mathcal{B}(\alpha)$, where $\mathbb{P}[A_\alpha(t)]$ is the probability that lattice site $\alpha$ is occupied at time $t$ and $\mathbb{P}[A_\alpha(t),A_\beta(t)]$ is the joint occupancy probability of neighbouring lattice sites $\alpha$ and $\beta$ at time $t$.

Mean-field models are widely used to predict ABM dynamics, however they fail to accurately predict dynamics in regions of parameter space in which the mean-field assumption is violated \citep{baker_correcting_2010,fadai_accurate_2019,hiebeler_stochastic_1997,matsiaka_continuum_2017,middleton_continuum_2014,newman_many-body_2004,simpson_distinguishing_2014}.  We depict ABM snapshots for two simulations of the BDM process in Figure \ref{fig:ABM_sims}. The mean-field assumption seems to be satisfied during the simulation with $(P_p,P_d)=(0.01,0.0025)$: agents appear uniformly distributed, which indicates that neighbouring site occupancies are independent of each other. As expected, we observe close agreement between $C(t)$ and $\mean{C_\text{ABM}(t)}$ for this parameter combination. For the simulation with $(P_p,P_d)=(1,0.25)$, however, the ABM simulation exhibits strong clustering. In this case, neighbouring site occupancies will be dependent within and outside of cluster regions, which violates the mean-field assumption. As a result, we observe poor agreement between the mean-field model and the ABM simulation output, $\mean{C_\text{ABM}(t)}$. Similar results have been documented previously for this model over a wide range of parameter values: the mean-field model matches ABM output well for small values of $P_p/P_m$ and $P_d/P_m$ ($P_m$ is fixed at unity in all simulations in this study), but this agreement worsens as either of these ratios increase \citep{baker_correcting_2010}.

\subsection{An SIR ABM}\label{sec:SIR_ABM}

Susceptible, infected, and recovered (SIR) models are used in epidemiology to model and predict the emergence of infectious  \citep{blackwood_introduction_2018} and waterborne \citep{eisenberg_identifiability_2013} diseases. In this section, we detail how the previous modelling framework can be extended to derive DE models for such ABMs of disease spread. We introduce the model rules in Section \ref{subsubsec:SIR_ABM_rules} and then derive the mean-field DE model and compare it to ABM output in Section \ref{subsubsec:SIR_CG}.

\subsubsection{ABM Rules}\label{subsubsec:SIR_ABM_rules}

 We use an equivalent lattice to that presented in Section \ref{sec:BMD_process}. Each lattice site, $\alpha$, can now take one of four states over time: $0_\alpha(t)$, $S_\alpha(t)$, $I_\alpha(t)$, and $R_\alpha(t)$ denote that $\alpha$ is unoccupied, or occupied by a \emph{susceptible}, \emph{infected}, or \emph{recovered} agent at time $t$, respectively. We assume three rules governing how agents move, infect, and recover in our SIR model. For agent movement, we assume that each agent moves with rate $P_m$. When an agent attempts to migrate, the agent chooses a neighbouring lattice site $\beta\in\mathcal{B}(\alpha)$ randomly to move to. If the chosen site is already occupied, then the migration event is aborted. We write this process as a bimolecular reaction with rate $P_m/4$:
\begin{equation}
    Y_\alpha + 0_\beta \xrightarrow{P_m/4} 0_\alpha + Y_\beta, \  \beta\in\mathcal{B}(\alpha), \label{eq:move_rate_SIR}
\end{equation}
for $Y\in \{S,I,R$\}.

The second rule governs infection of agents, which occurs with rate $P_I$. During an infection event, an infected agent at lattice site $\alpha$ will randomly infect an agent at a neighbouring lattice site $\beta\in\mathcal{B}(\alpha)$. If the chosen site $\beta$ is occupied by a susceptible agent, then the susceptible agent becomes infected. Otherwise, the infection does not alter the state of lattice site $\beta$. We model this rule using a bimolecular reaction with rate $P_I/4$:
\begin{equation}
    I_\alpha + S_\beta \xrightarrow{P_I/4} I_\alpha + I_\beta, \  \beta\in\mathcal{B}(\alpha).\label{eq:infect_rate}
\end{equation}
The final rule concerns the recovery of infected agents: infected agent recover with rate $P_R$. We model this process using a monomolecular reaction with rate $P_R$:
\begin{equation}
    I_\alpha \xrightarrow{P_R} R_\alpha. \label{eq:recover_rate}
\end{equation}

Simulation of the SIR ABM proceeds as follows. We begin each simulation by randomly placing susceptible agents in 49\% of the lattice sites, infected agents in 1\% of the lattice sites, and leaving the remaining lattice sites unoccupied. \new{Note that each ABM simulation begins with a new initial configuration of agent locations throughout the lattice}. Because there is no death or birth in the model, the proportion of occupied lattice sites is fixed at $M=0.5$ for all time. We use reflecting boundary conditions, which model a no-flux condition in the spatial domain. The Gillespie Algorithm is used to simulate the model \citep{gillespie_exact_1977}. From the $n^\text{th}$ of $N$ identically prepared simulations, we let $S^{(n)}_\text{ABM}(t)$, $I^{(n)}_\text{ABM}(t)$, and $R^{(n)}_\text{ABM}(t)$, denote the fractions of susceptible, infected, and recovered agents in the model over time, respectively (e.g., $S^{(n)}_\text{ABM}(t)$ is equal to the number of susceptible agents at time $t$ divided by $MX^2$). We then estimate the averaged ABM fraction for each subpopulation by averaging over all $N$ simulations
\begin{equation}
    \mean{S_\text{ABM}(t)} = \dfrac{1}{N}\sum_{n=1}^N S_\text{ABM}^{(n)}(t); \ \ \ \ \ 
    \mean{I_\text{ABM}(t)} = \dfrac{1}{N}\sum_{n=1}^N I_\text{ABM}^{(n)}(t); \ \ \ \ \
    \mean{R_\text{ABM}(t)} = \dfrac{1}{N}\sum_{n=1}^N R_\text{ABM}^{(n)}(t).
\end{equation}
We depict snapshots of two simulations of the SIR model in Figure \ref{fig:SIR_ABM_compare}. In both cases, we observe that the small initial proportion of infected agents causes an outbreak of infection. The majority of agents have become infected and then recovered by the end of both simulations.

\subsubsection{Model coarse graining} \label{subsubsec:SIR_CG}
 
\begin{figure}
    \centering
    \includegraphics[width=\textwidth]{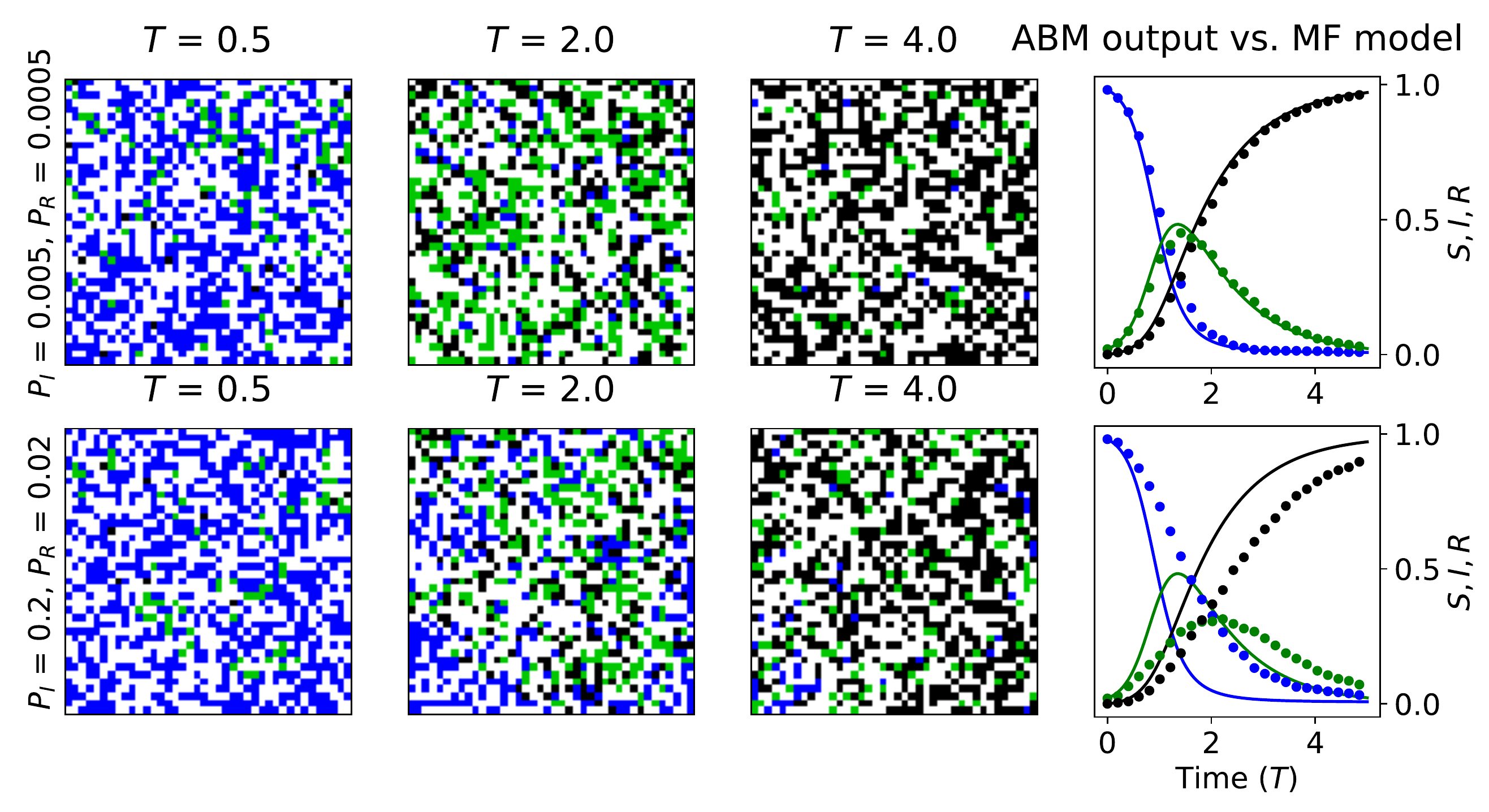}
    \caption{Simulation snap shots for the SIR ABM. Blue pixels denote $S_\alpha(t)$, green pixels denote $I_\alpha(t)$, black pixels denote $R_\alpha(t)$, and white pixels denote $0_\alpha(t)$. In the rightmost column, we compare predictions of the mean-field SIR model (Equation \eqref{eq:SIR_MF}) to the \new{computed ratios of $S,I,\text{ and }R$ from one ABM simulation}. Solid lines correspond to the mean-field (MF) model and dots correspond to ABM simulation output. Quantities are plotted against nondimensionalised time $T=P_Rt$ for ease of interpretation. The ABM was computed on a square lattice of length $X=40$.}
    \label{fig:SIR_ABM_compare}
\end{figure}

 We show in Appendix \ref{sec:app_deriving_DE_SIR} that the mean-field model for the SIR process is given by the frequently-used system of equations:
 \begin{equation}
    \dfrac{\text{d}S}{\text{d}t} = -M P_ISI; \ \ \ \ \  \dfrac{\text{d}I}{\text{d}t} = M P_ISI - P_RI; \ \ \ \ \  \dfrac{\text{d}R}{\text{d}t} = P_RI.\label{eq:SIR_MF}
\end{equation}
In Equation \eqref{eq:SIR_MF}, the variables $S(t),I(t),R(t)$ denote the fraction of susceptible, infected, and recovered agents at time $t$, respectively. In Figure \ref{fig:SIR_ABM_compare}, we depict snapshots from two simulations of the SIR ABM together with evolution of the corresponding ABM densities. The ABM simulation with $(P_I,P_R)$ = (0.005,0.0005) appears well-mixed for all agent-types, which would satisfy the mean-field assumption. As a result, simulations of the mean-field model predict the ABM density well. If we increase the infection and recovery rates to $(P_I,P_R)$ = (0.2,0.02), however, then the resulting ABM simulation has separate patches comprised primarily of infected agents or susceptible agents.  These patches of single agent types decrease the population-wide average infection rate because infected agents in the middle of an infected cluster are unable to infect any susceptibles. As a result, the mean-field assumption is violated within these patches, and the mean-field model cannot accurately predict ABM dynamics.   

For both of the BDM and SIR models, we have seen that some parameter regimes lead to close agreement between the ABM output and mean-field model predictions, whereas other parameters lead to poor agreement between the two. There is thus a need to develop methods that can determine when mean-field modes will accurately predict ABM dynamics and find novel DE models that accurately predict ABM dynamics when the mean-field models fail to do so.

\section{Equation learning} \label{sec:EQL}

For many EQL studies, the goal is to infer a dynamical systems model, written as $C(t)$, from a time-varying dataset, $C_d(t)$. This dynamical system can broadly be written as
\begin{equation}\label{eq:dynamical_system}
    \dfrac{\text{d}C(t)}{\text{d}t} = \mathcal{F},
\end{equation}
where  $\mathcal{F}$ describes the dynamics of $C(t)$. When $C_d(t)$ is a time-varying scalar quantity (or a vector of scalar quantities), then an ODE model is relevant, in which case $\mathcal{F}=\mathcal{F}(t,C)$. When $C_d(t)$ varies over time and a one-dimensional spatial dimension, $x$, then a PDE model may be more relevant, in which case $\mathcal{F}=\mathcal{F}(t,x,C,C_x,C_{xx},\dots)$. In the following sections, we exemplify how methods from EQL may be used to learn a form of $\mathcal{F}$ from output ABM data. 

\begin{figure}
    \centering
    \includegraphics[width=0.99\textwidth]{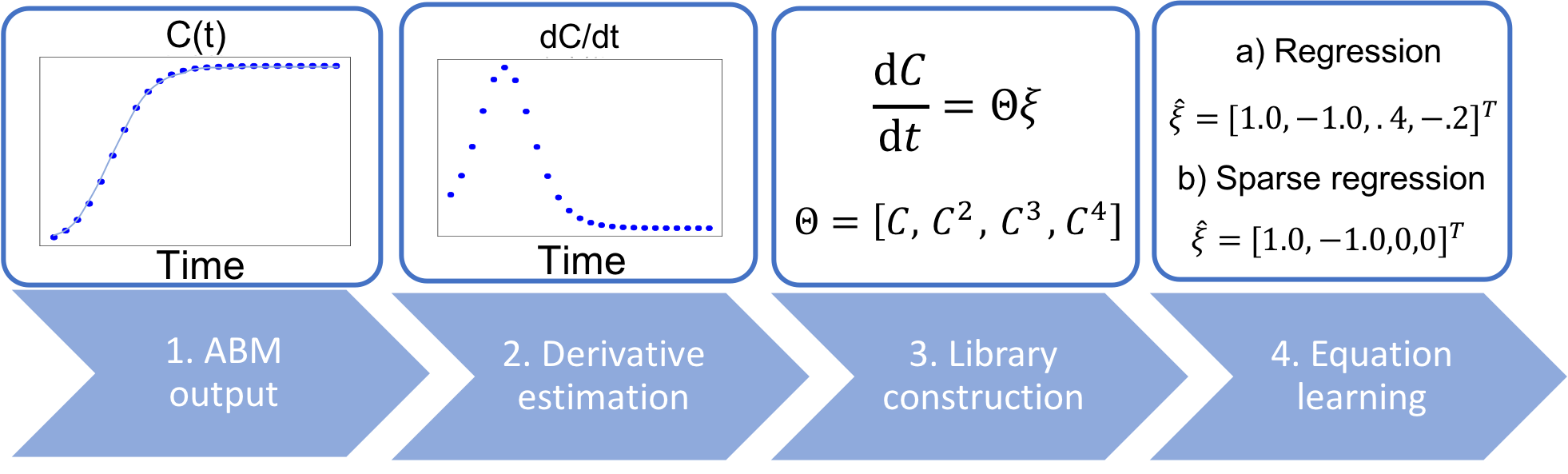}
    \caption{EQL pipeline. Step 1: Generate averaged ABM output; Step 2: Estimate the temporal derivative of the ABM output; Step 3: Library construction; Step 4: Equation inference. At Step 4, one can either perform (a) regression  or (b) sparse linear regression  to learn an equation for the ABM output. We will consider both the Lasso and Greedy sparse regression algorithms to perform EQL.}
    \label{fig:EQL_ABM_tutorial}
\end{figure}

\subsection{Model learning example}\label{subsec:ABM_DE_case_study}

In this section, we outline the steps one may take to learn a DE model from ABM data for the BDM model (Figure \ref{fig:EQL_ABM_tutorial}). Code is provided for the results presented in this section in the file {\color{blue} EQL Tutorial.ipynb}. 

\subsubsection{EQL pipeline}\label{subsubsec:EQL_pipeline}

We  illustrate how to use EQL methods using data from the BDM process described in Section \ref{sec:BMD_process}. In the first step, we simulate the ABM 50 times with parameter values ($P_p,P_m,P_d$)=(0.01,0.005,1.0) and average the output to acquire $C_d(t)=\mean{C_\text{ABM}(t)}$. The time vector, $t$, is sampled on an equispaced grid such that $t_i = (i-1)\Delta t, \ i = 1,\dots,100$ for some small $\Delta t > 0$.

In the second step of this process, we estimate the numerical derivative of $C_d(t)$. Finite difference computations are a simple method to approximate derivatives \citep{leveque_finite_2007}. We use centered differences at the internal time points and forward and backward differences at the first and final time points, respectively:
\begin{align}
     \dfrac{\text{d}C_d(t_1)}{\text{d}t} &= \dfrac{C_d(t_{2}) - C_d(t_{1})}{\Delta t} , \nonumber \\
     \dfrac{\text{d}C_d(t_i)}{\text{d}t} &= \dfrac{C_d(t_{i+1}) - C_d(t_{i-1})}{2\Delta t},\ \ \  i = 2,\dots,n-1,  \nonumber\\
     \dfrac{\text{d}C_d(t_n)}{\text{d}t} &= \dfrac{C_d(t_n) - C_d(t_{n-1})}{\Delta t}. 
 \end{align}
The resulting computation is plotted in Step 2 of the EQL pipeline depicted in Figure \ref{fig:EQL_ABM_tutorial}. 

The third step of the EQL pipeline requires the construction of a library of potential terms for inclusion in the inferred DE model. We saw that polynomials in $C(t)$ can describe the ABM output in Section \ref{sec:ex1_logistic_deriv}, so we assume that the underlying model here is a polynomial in $C(t)$. Recall from the rules of the BDM ABM  that each agent interacts with its four neighbouring sites, so we further assume that this polynomial is up to fourth order. \new{See reference \citep{simpson_cell_2010} for scenarios where fourth or higher order polynomials are needed to match ABM output for reactions involving two agents.} As a nonzero constant in $\mathcal{F}$ would represent a constant source or sink of agents (which is not present in the ABM), we set the constant in this polynomial to be zero. Altogether, we propose the following possible model for the BDM process
\begin{equation}
    \dfrac{\text{d}C}{\text{d}t} = \sum_{i=1}^4 \xi_i C^i,\label{eq:EQL_proposed_model}
\end{equation}
for unknowns $\xi_i\in\mathbb{R}$. Given data $C_d(t)$, we substitute into Equation \eqref{eq:EQL_proposed_model} to arrive at the following linear system of equations satisfied by the unknowns $\xi_1,\dots,\xi_4$:

\begin{equation}
    \begin{bmatrix}
    \dfrac{\text{d}C_d(t_1)}{\text{d}t} \\ \dfrac{\text{d}C_d(t_2)}{\text{d}t} \\  \vdots \\ \dfrac{\text{d}C_d(t_n)}{\text{d}t} \\ \end{bmatrix}
    =  \begin{bmatrix}
    C_d(t_1) & C^2_d(t_1) & C^3_d(t_1) & C^4_d(t_1) \\
    C_d(t_2) & C^2_d(t_2) & C^3_d(t_2) & C^4_d(t_2) \\
    \vdots & \vdots & \vdots & \vdots \\
    C_d(t_n) & C^2_d(t_n) & C^3_d(t_n) & C^4_d(t_n) \\
    \end{bmatrix}
    \begin{bmatrix}
    \xi_1 \\ \xi_2 \\  \xi_3 \\ \xi_4 \\ \end{bmatrix}.\label{eq:ABM_system_long}
\end{equation}
For convenience, we re-write Equation \eqref{eq:ABM_system_long} as 

\begin{equation}
    \dfrac{\text{d}C_d}{\text{d}t} = \bm{\Theta}\xi, \label{eq:ABM_system}
\end{equation}
where the columns of the $n\times4$ matrix $\bm{\Theta}$ contain the library terms evaluated at the data points $C_d(t_i),i=1,\dots,n$.

\begin{figure}
    \centering
    \includegraphics[width=0.45\textwidth]{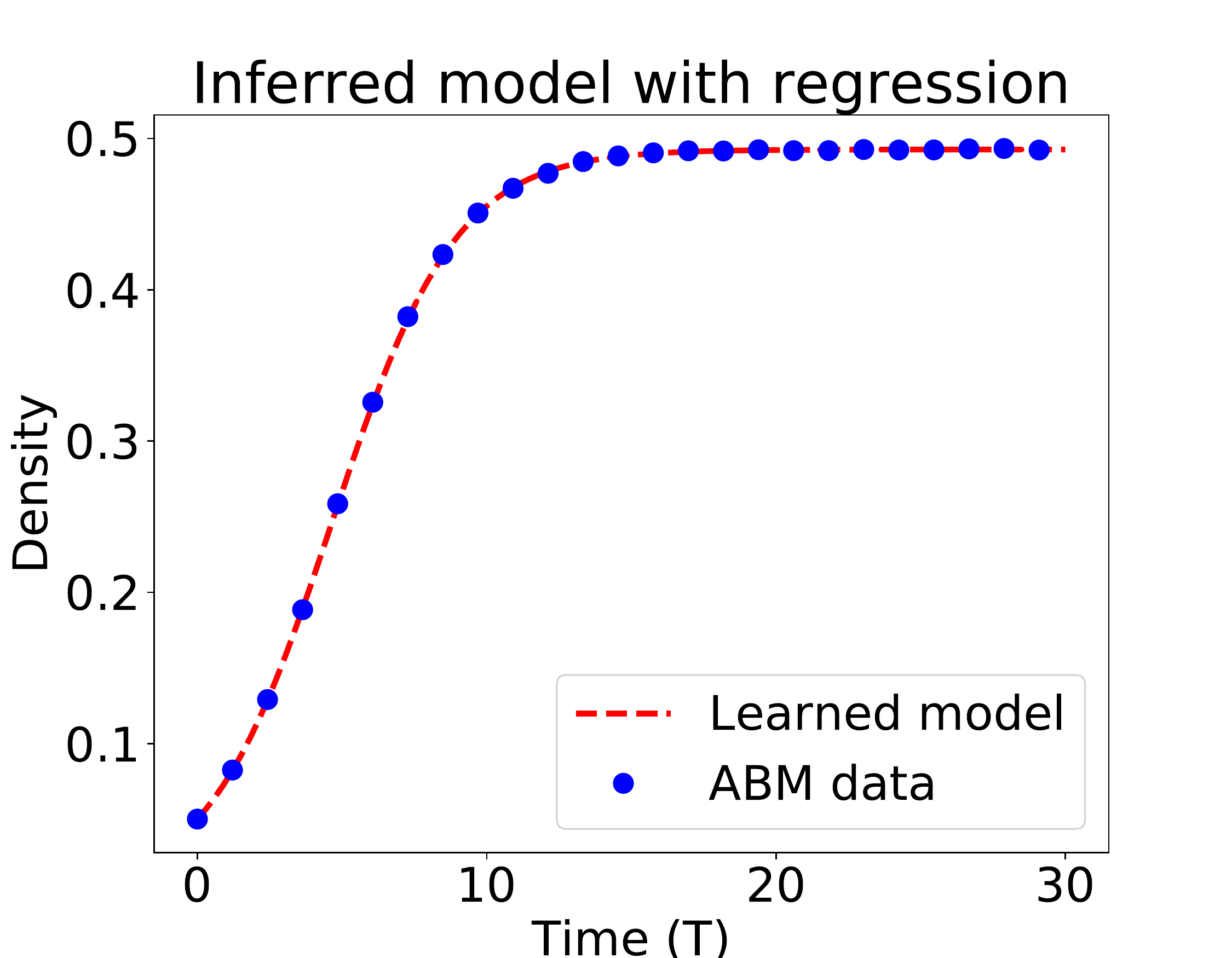}
    \includegraphics[width=0.45\textwidth]{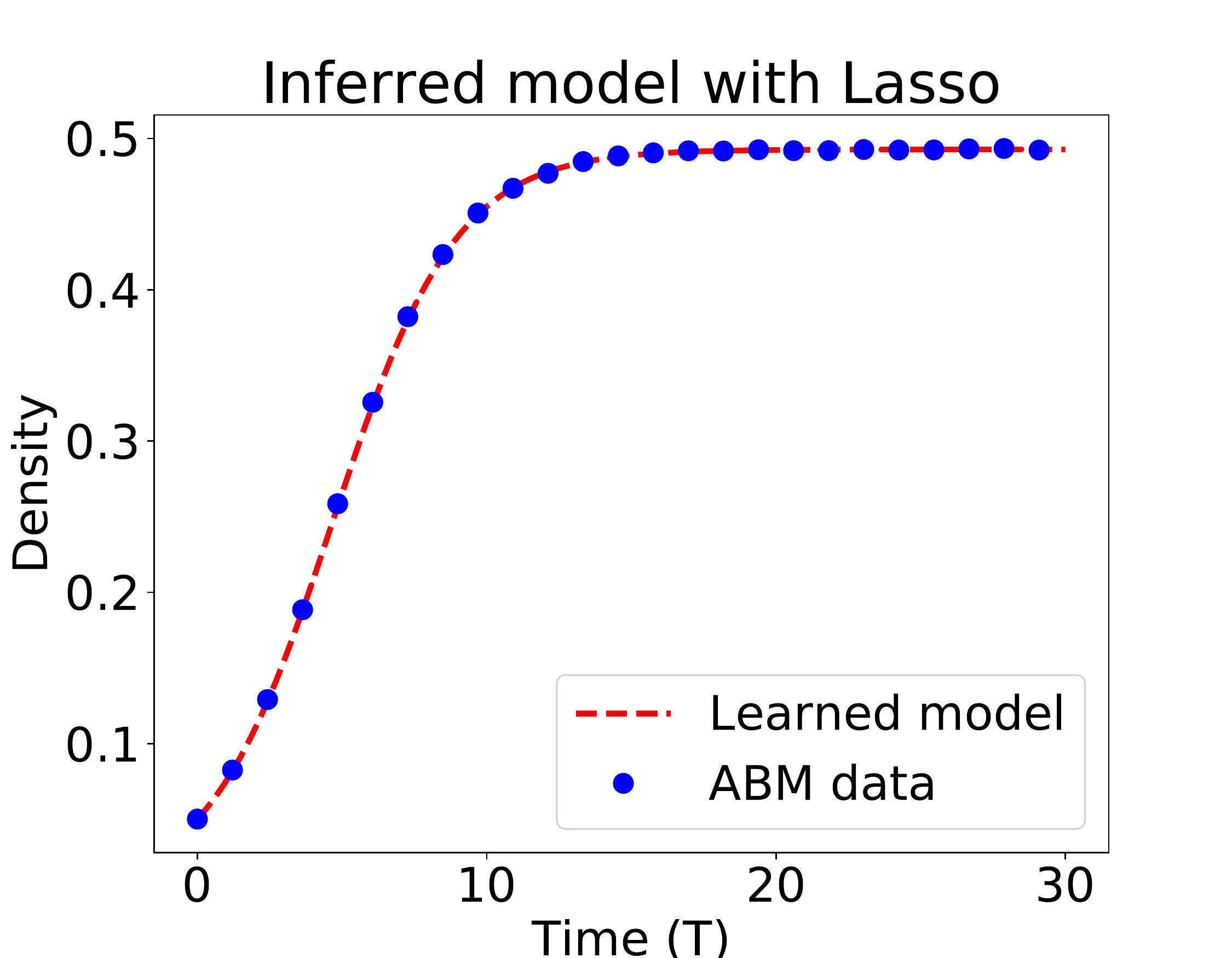}
    \caption{Simulations of the inferred models for the BDM ABM data. (left) The inferred DE model $\text{d}C/\text{d}t=0.0100C - 0.0182C^2 + 0.0037C^3 - 0.0019C^4$  (determined by solving Equation \eqref{eq:ABM_system} using linear regression) is depicted against output ABM data, $C_d(t)$. (Right) The inferred DE model $\text{d}C/\text{d}t=0.0047C - 0.0095C^2$ (determined by solving Equation \eqref{eq:ABM_system} using Lasso) is depicted against $C_d(t)$.}
    \label{fig:tutorial_sim_models}
\end{figure}

The fourth step in the EQL pipeline is to infer the form of the DE model that best approximates the dynamics of $C_d(t)$. We do so by finding the least-squares solution of Equation \eqref{eq:ABM_system}, given by
\begin{equation}
    \hat{\xi} = \arg\min_{\xi\in\mathbb{R}^4}\left\{ \dfrac{1}{n}\left\|\dfrac{\text{d}C(t)}{\text{d}t} - \bm{\Theta}\xi\right\|_2^2\right\} \label{eq:regression}.
\end{equation}
 We solve Equation \eqref{eq:regression} using numpy's \textbf{lstsq} command from the linear algebra package and find $\xi = [0.0048,-0.0105,0.0031,-0.0030]^T$. This solution suggests that the following model best describes the ABM dynamics:
\begin{equation}
    \dfrac{\text{d}C(t)}{\text{d}t} = 0.0048C-0.0105C^2 + 0.0031C^3 -0.0030C^4.\label{eq:inferred_ABM_EQ_full}
\end{equation}
We numerically simulate Equation \eqref{eq:inferred_ABM_EQ_full} with initial condition $C(0)=C_d(0)$ using a fourth order Runge-Kutta method \citep{leveque_finite_2007}. The resulting output, $C(t)$, is depicted against against $C_d(t)$ in Figure \ref{fig:tutorial_sim_models} and we observe that this model accurately predicts the ABM output.

By solving Equation \eqref{eq:ABM_system} directly, we are likely to find forms of $\mathcal{F}$ with many terms because the system is overdetermined when $n$, the number of data points, satisfies $n\gg4$. We may wonder if a simpler form of $\mathcal{F}$ can also accurately describe the ABM dynamics. A logistic model such as that presented in Equation \eqref{eq:logistic} has only two terms and may describe the data accurately, for example. Instead of solving Equation \eqref{eq:regression}, we can use sparse regression methods to find a sparse vector, $\xi$, to solve Equation \eqref{eq:ABM_system}. There are many approaches to sparse regression, including the \emph{least absolute shrinkage and selection operator} (Lasso), ridge regression, and the Greedy algorithm \citep{tibshirani_regression_1996,rudy_data-driven_2017,zhang2009adaptive}. We utilize the Lasso algorithm in this section, but will return later to a discussion of whether alternative methods should also be considered. The Lasso method solves the regularized system
\begin{equation}
    \hat{\xi} = \arg\min_{\xi\in\mathbb{R}^4}\left\{ \dfrac{1}{n}\left\|\dfrac{\text{d}C(t)}{\text{d}t} - \bm{\Theta}\xi\right\|_2^2 + \lambda\|\xi\|_1 \right\}, \label{eq:lasso}
\end{equation}
for some $\lambda>0$, which is called a \emph{regularization parameter} \citep{tibshirani_regression_1996}. Regularization is used to avoid extreme values in $\xi$ and to prevent overfitting. There are many approaches to solve the Lasso problem; here we use the  Fast iterative Shrinking-Thresholding Algorithm (FISTA) \citep{beck_fast_2009}. We provide the pseudo-code for this algorithm in  Algorithm \ref{algo:lasso} of Appendix \ref{sec:fista}. Using a regularization strength of $\lambda=0.0004$ (see Appendix \ref{sec:hyperselectlasso} for a discussion of hyperparameter selection) we find the resulting vector to be $\hat{\xi}=[0.0047,-0.0095,0,0]$. This estimate suggests that the model equation
\begin{equation}
    \dfrac{\text{d}C}{\text{d}t} = 0.0047C-0.0095C^2, \label{eq:learned_eqn_lasso}
\end{equation}
should be a more parsimonious model for the BDM process than Equation \eqref{eq:inferred_ABM_EQ_full}. We depict the solution to Equation \eqref{eq:learned_eqn_lasso} (with initial condition $C(0)=C_d(0)$) against the ABM output in Figure \ref{fig:tutorial_sim_models} and observe that this DE model accurately approximates $C_d(t)$. Furthermore, we observe that the form of Equation \eqref{eq:learned_eqn_lasso} is similar to the logistic DE, 
\begin{equation}
    \dfrac{\text{d}C}{\text{d}t} = P_pC(1-C) - P_dC = (P_p-P_d)C-P_pC^2 .\label{eq:logistic_restate}
\end{equation}
By comparing coefficients between Equations \eqref{eq:learned_eqn_lasso} and \eqref{eq:logistic_restate} we can estimate the mechanistic ABM parameters $P_p$ and $P_d$ as $\hat{P_p}=0.0095\text{ and }\hat{P_d}=0.0048$. These estimates are very close to the true underlying values of $P_p=0.01\text{ and }P_d=0.05$. The proposed EQL methodology is thus able to simultaneously infer a DE model that accurately predicts ABM output and provides realistic parameter estimates when combined with the mean-field model \new{(we note that more common forms of parameter estimation, such as maximum likelihood, may need to be used after the equation form has been determined)}. We have thus shown that concepts from EQL can be used to determine simple DE models that accurately describe ABM dynamics.

\subsubsection{Different forms of sparse regression}

As mentioned previously, there are many approaches to find sparse solutions to Equation \eqref{eq:ABM_system}. The Greedy algorithm is another popular algorithm for sparse regression and solves a similar problem to the Lasso problem from Equation \eqref{eq:lasso}. In the Greedy algorithm, we solve
\begin{equation}
    \hat{\xi} = \arg\min_{\xi\in\mathbb{R}^4}\left\{ \dfrac{1}{n}\left\|\dfrac{\text{d}C(t)}{\text{d}t} - \bm{\Theta}\xi\right\|_2^2 + \lambda\|\xi\|_0 \right\}, \label{eq:greedy}
\end{equation}
for some regularization parameter $\lambda>0$. Here, $\|\xi\|_0$ counts the number of nonzero terms in $\xi$. We use the forward-backward approach to solve this system, which converts the regularization parameter $\lambda$ into a tolerance hyperparameter. Pseudo-code for this algorithm is provided in \citep{zhang2009adaptive}. We use this algorithm on the ABM data with a tolerance of $0.0001$ and find the resulting vector to be $\hat{\xi}=[0.0047,-0.0095,0.0001,0]$, which suggests that the model equation
\begin{equation}
    \dfrac{\text{d}C}{\text{d}t} = 0.0047C-0.0095C^2 + 0.0001C^3, \label{eq:learned_eqn_greedy}
\end{equation}
is able to describe the ABM dynamics. We note that this learned equation is similar in form to the learned equation from Lasso (Equation \eqref{eq:learned_eqn_lasso}), although there is an extra cubic term with a small coefficient. We will consider both the Lasso and Greedy algorithms for EQL in future case studies.

\section{Case studies for EQL in ABM analysis}\label{sec:case_study}

We use this section to explore how methods from EQL can aid modellers in performing ABM analysis with EQL methods. We will do so through five case studies pertaining to: how learned equations change with ABM parameters in Section \ref{subsec:case_study_logistic}, how learned equations are affected by the number of performed ABM simulations in Section \ref{subsec:case_study_2_reals}, the performance of learned equations in predicting unobserved ABM dynamics in Section \ref{subsec:case_study_3_unobserved}, the performance of EQL methods for DE model selection from ABM data in Section \ref{subsec:case_study_4_model_selection}, and the performance of EQL methods for learning systems of equations in Section \ref{subsec:varying_infection_SIR}. Through these case studies, we introduce and address six questions on the use and efficacy of EQL for ABM analysis.

\subsection{Case study 1: Comparing DE models in describing ABM dynamics} \label{subsec:case_study_logistic}

Coarse-grained models are advantageous for ABM analysis because they are easy to interpret, formulate, and solve. Unfortunately, coarse-grained DE models only provide accurate ABM analysis in parameter regimes where ABM simulations adhere to their underlying assumptions \citep{fadai_accurate_2019,simpson_distinguishing_2014}. Previous studies have defined criteria to aid modellers in determining when to rely on mean-field models, \new{but these approaches are often only valid for simple ABMs and heuristic in nature for more complex scenarios, such as bistable systems} \old{such as those detailed in this work, but such criteria are often heuristic in nature} \citep{baker_correcting_2010,fadai_accurate_2019,grima_effective_2010,simpson_distinguishing_2014,smith_analytical_2016}. The purpose of this case study is to determine if EQL methods can be used as a simple test to determine when mean-field models accurately predict ABM dynamics, and to propose novel models for more accurate inference. These goals are summarized in the following questions
\begin{displayquote}
(Q1) Can EQL aid in determining when mean-field models accurately approximate ABM dynamics?\\
(Q2) Can methods from EQL discover novel DE models for accurate ABM analysis when the mean-field assumption is invalid?
\end{displayquote}
To address both Questions (Q1) and (Q2), we will test the ability of the mean-field and learned DE models to predict ABM dynamics for the BDM model over a range of mechanistic ABM parameter values. The SIR model is considered in a later case study.

\subsubsection{Varying proliferation rates for the BDM process}

We now investigate the performance of the mean-field and learned DE models in describing ABM output from the BDM process over a range of proliferation rates. The proliferation rates are varied as $P_p=0.01,0.05,0.1,0.5$. For each value of $P_p$, we set the death rate, $P_d$, to be half of $P_p$ and fix the migration rate at $P_m=1$. The ABM output is comprised of $C_d(t) = \mean{C_\text{ABM}(t)}$, which is averaged over $N=50$ ABM simulations to ensure convergence to mean behaviour and reduce the impact of noise. For model learning, we use the library of right-hand side terms $\bm{\Theta}=[C,\ C^2,\ C^3,\ C^4]$, and use the Greedy algorithm \citep{zhang2009adaptive} to sparsely solve the linear system $\text{d}C_d/\text{d}t=\bm{\Theta}\xi$. The code for this case study is provided in the file {\color{blue} \textbf{Case study 1 Varying parameters for BDM.ipynb}}.

\begin{figure}
    \centering
    \includegraphics[width=0.45\textwidth]{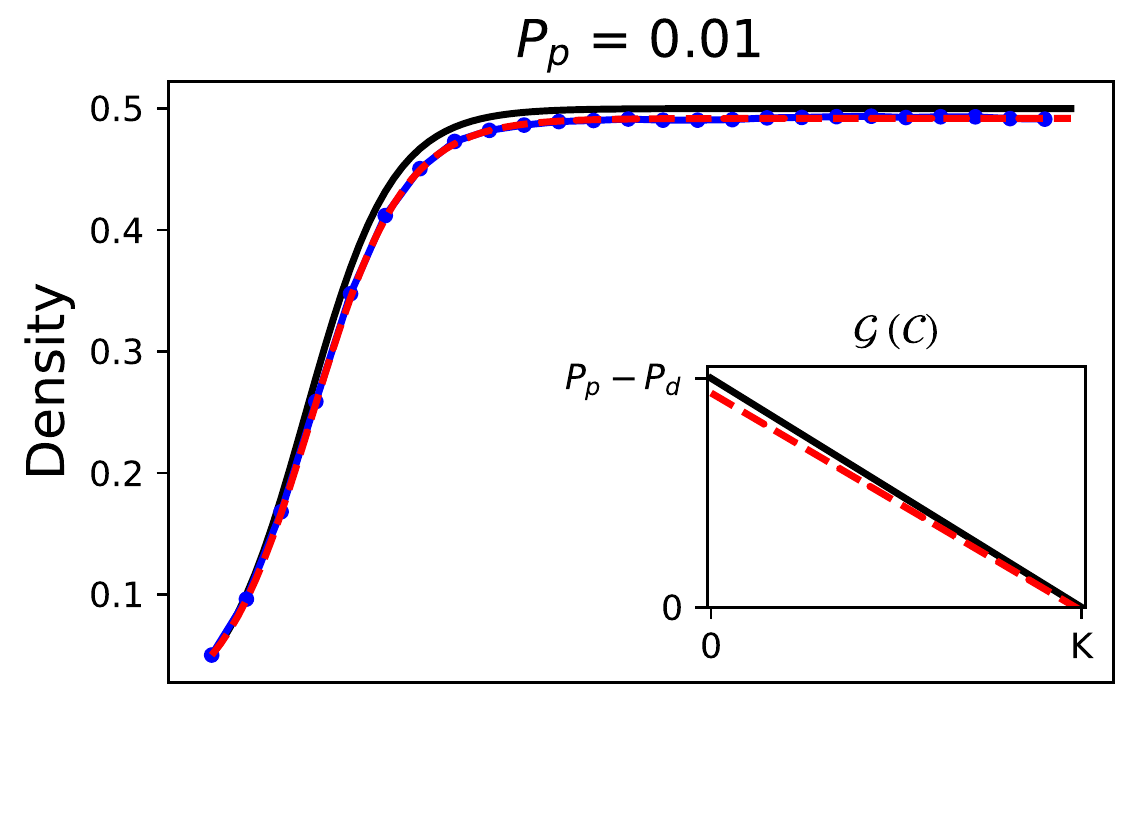}
    \includegraphics[width=0.45\textwidth]{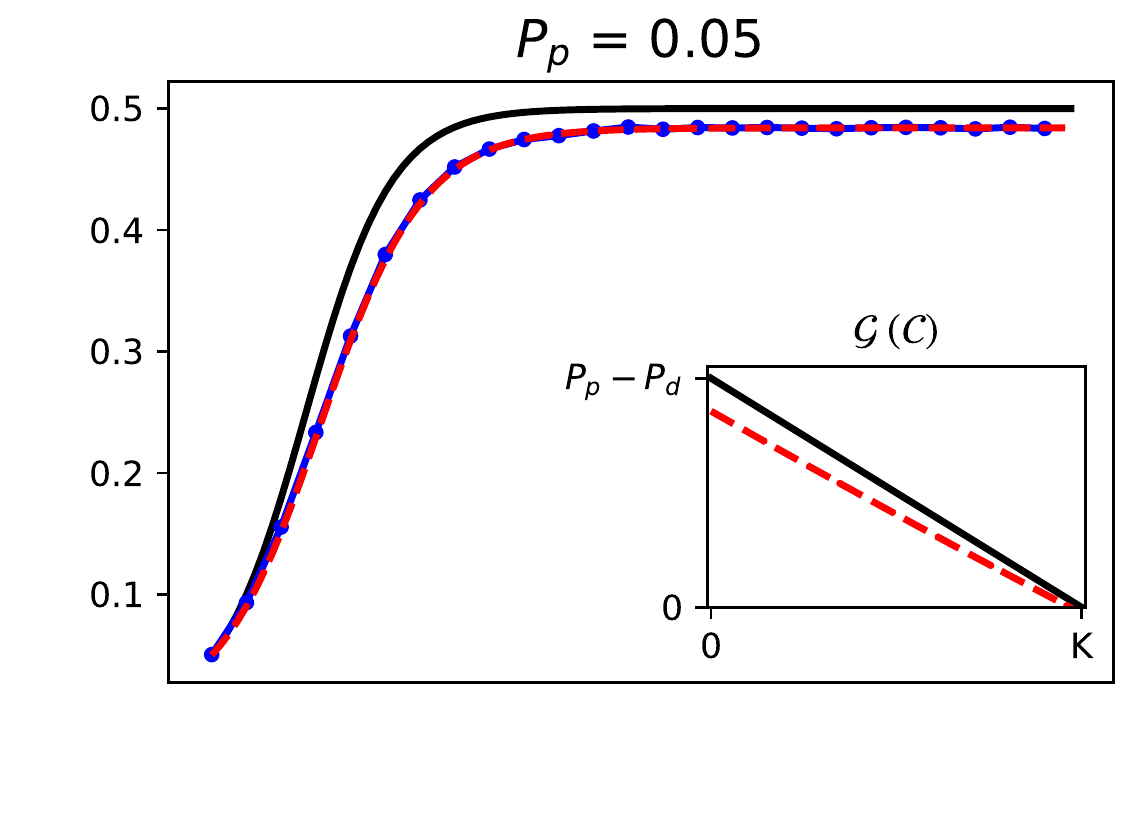}
    \includegraphics[width=0.45\textwidth]{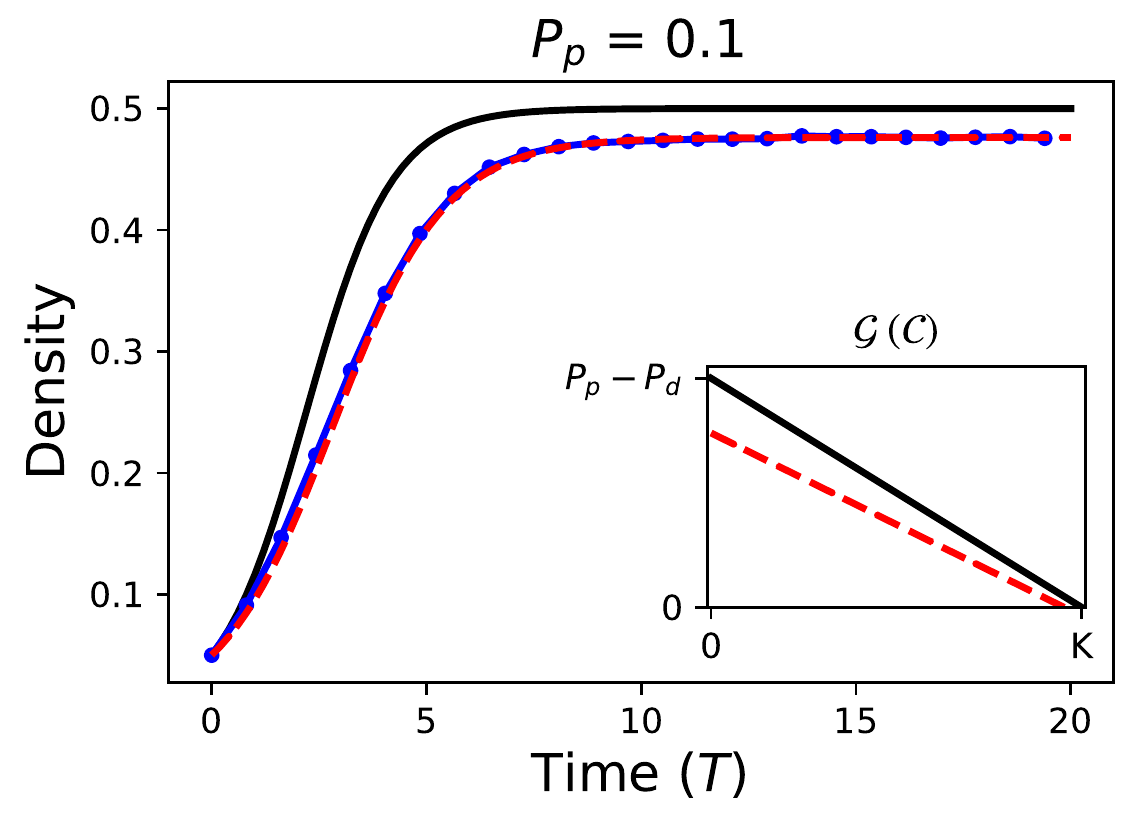}
    \includegraphics[width=0.45\textwidth]{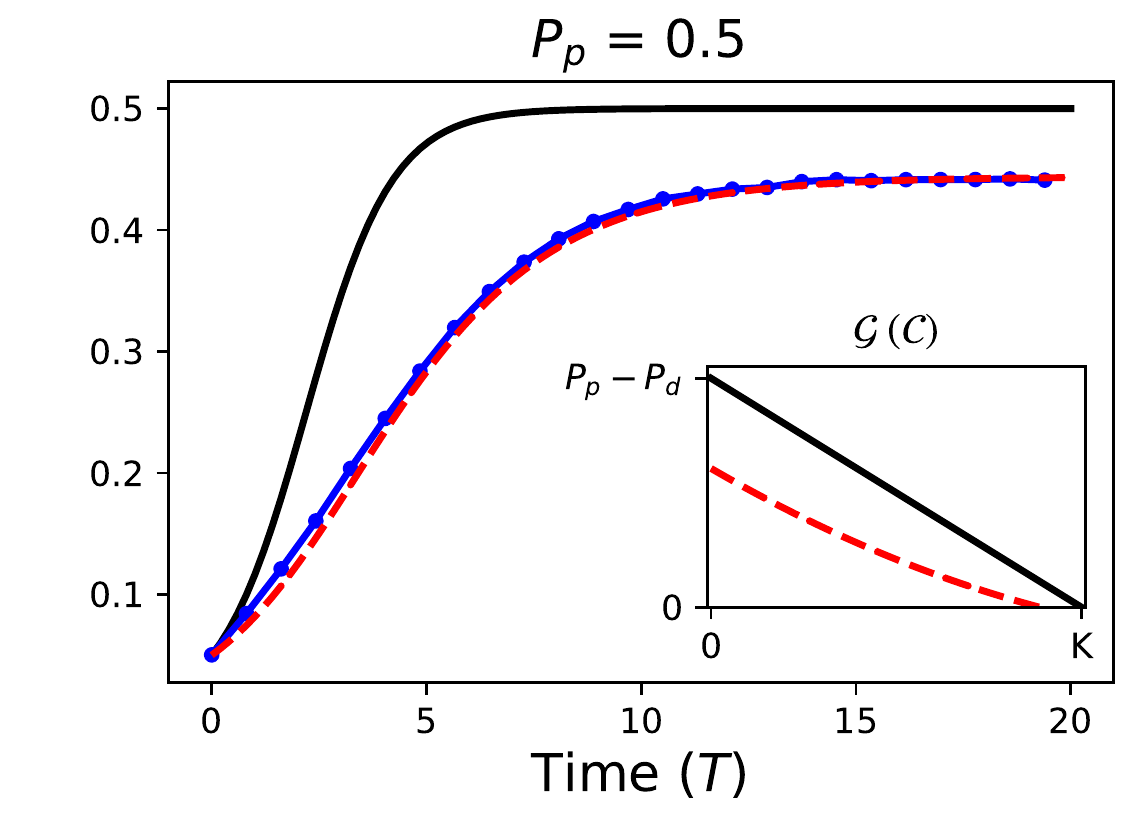}
    \caption{Case study 1. Comparing mean-field and learned model predictions to ABM data from the BDM process. In each figure we depict $\mean{C_\text{ABM}(t)}$ (blue dots and line) against the corresponding mean-field model (solid black line) and learned equation (red dashed line). All simulations are depicted as a function of nondimensionalised time $T=t(P_p-P_d)$. The insets in all frames depict the predicted per-capita growth rates, $\mathcal{G}(C)$, from both models where $\text{d}C/\text{d}t=C\mathcal{G}(C)$.}
    \label{fig:CS_1_BMD_data_comparison}
\end{figure}

We compute predictions of the mean-field and learned DE models with $\mean{C_\text{ABM}(t)}$ in Figure \ref{fig:CS_1_BMD_data_comparison} as well as the model equations and their mean-squared error (MSE) in approximating $\mean{C_\text{ABM}(t)}$ in Table \ref{tab:CS_table_1}. Both the mean-field model and learned DE model provide accurate predictions of the ABM output for $P_p=0.01$, with a MSE between the simulated model and $\mean{C_\text{ABM}(t)}$ of 0.0011 and 0.0001, respectively. As $P_p$ increases to 0.05, 0.1, and 0.5, the mean-field model does an increasingly poor job in predicting the ABM data, overpredicting agent density for all time. The MSE of the mean-field DE model increases with $P_p$. The learned DE models, on the other hand, accurately predict the ABM data and maintain MSE values below 0.0005 for all values of $P_p$. The learned DE model form is similar to the mean-field model for $P_p=0.01$, but for $P_p=0.05,0.1,0.5$ the learned model also recovers cubic terms. When the learned model resembles the mean-field model, then the mean-field model accurately predicts $\mean{C_\text{ABM}(t)}$. On the other hand, when the learned model deviates from the mean-field model, the mean-field model poorly predicts the ABM data. 

We suggest that the mean-field model can make accurate predictions of ABM behaviours when the learned equation closely resembles the mean-field model \new{(including both equation form and parameter estimates)}; otherwise, the mean-field model can lead to inaccurate  predictions of ABM behaviours.  Consider the \emph{per-capita growth rate} as one illustrative example. For a DE model of the form $\text{d}C/dt=\mathcal{F}(C)$, the per capita growth rate is defined by $\mathcal{F}(C) = C\mathcal{G}(C)$, and it quantifies the average contribution of each individual to population growth over time. We plot $\mathcal{G}(C)$ for the mean-field model and each learned model in the insets of Figure \ref{fig:CS_1_BMD_data_comparison}. The mean-field model predicts that the per-capita growth is a linear decreasing function connecting (0,$P_p-P_d$) and ($K$,0), where $K=(P_p-P_d)/P_p$ is the carrying capacity predicted by the mean-field model. The learned model predictions of $\mathcal{G}(C)$ closely resemble the mean-field model predictions of $\mathcal{G}(C)$ for $P_p=0.01$ and 0.05. At the larger proliferation rates of $P_p=0.1$ and 0.5, however, the learned model per-capita growth rates are much lower than the mean-field model rates. Recall that higher rates of proliferation lead to spatial clustering of agents in the ABM: this clustering reduces the averaged per-capita growth rate, which the learned model can account for but the mean-field model does not. The effective carrying capacity of the ABM reduces as $P_p$ increases (with $P_d=P_p/2$), which is again likely due to increased spatial clustering. All learned models accurately capture this reduction in the carrying capacity, whereas the mean-field models do not. 

The EQL pipeline results for these four simulated datasets suggests that the mean-field model will accurately describe ABM data when the learned equation form matches that of the mean-field model. When the mean-field models are not able to accurately predict ABM dynamics, learned models of the form
\begin{equation}
    \dfrac{\text{d} C}{\text{d}t} = \alpha C +  \beta C^2 + \gamma C^3, \ \ \  \alpha,\beta,\gamma \in\mathbb{R}, \label{eq:learned_CS1}
\end{equation}
can accurately predict ABM data instead. This form of learned equation was able to provide more accurate ABM analysis than the mean-field model over a range of parameter values.

\begin{table}[]
    \centering
    \begin{tabular}{|c|c|l|l|}
    \hline
    $P_p$ & \ $P_d$ & Mean-field Model (MSE) & Learned Model (MSE) \\ 
    \hline
    $0.01$ & $0.005$ & $\nicefrac{\text{d}C}{\text{d}t} = 0.005C    - 0.01C^2   $ (0.0011) & $\nicefrac{\text{d}C}{\text{d}t} = 0.00468C    - 0.0095C^2    - 0.0C^3   $ (0.0001) \\
    \hline
    $0.05$ & $0.025$ & $\nicefrac{\text{d}C}{\text{d}t} = 0.025C    - 0.05C^2   $ (0.0026) & $\nicefrac{\text{d}C}{\text{d}t} = 0.02134C    - 0.04572C^2    + 0.00343C^3   $ (0.0002) \\
    \hline
    $0.1$ & $0.05$ & $\nicefrac{\text{d}C}{\text{d}t} = 0.05C    - 0.1C^2   $ (0.004) & $\nicefrac{\text{d}C}{\text{d}t} = 0.03962C    - 0.09057C^2    + 0.01561C^3   $ (0.0003) \\
    \hline
    $0.5$ & $0.25$ & $\nicefrac{\text{d}C}{\text{d}t} = 0.25C    - 0.5C^2   $ (0.01) & $\nicefrac{\text{d}C}{\text{d}t} = 0.15671C    - 0.49984C^2    + 0.33125C^3   $ (0.0005) \\
    \hline
    \end{tabular}
    \caption{Case study 1. Mean-field and learned DE models for the BDM process for various values of $(P_p,P_d)$. MSE denotes the mean squared error between the model solution and $\mean{C_\text{ABM}(t)}$.}
    \label{tab:CS_table_1}
\end{table}

\subsection{Case study 2: Altering the number of ABM simulations} \label{subsec:case_study_2_reals}

ABMs are inherently noisy due to the random updating of agent states that occurs during simulation. When using EQL methods to analyse such ABMs, one should take care to ensure they learn the mean dynamics and do not overfit to small trends in the data. We averaged ABM data over a large number of ABM simulations in previous sections of this study to ensure the data had converged to its mean value. Performing such extensive simulations may not be feasible for computationally intensive ABMs, however, so we now investigate how the learned DE model changes with the number of ABM simulations. This can be summarized with the following question:
\begin{displayquote}
(Q3) How can we determine when enough ABM simulations have been performed for accurate DE model learning?
\end{displayquote}

To investigate (Q3), we consider data sets from the BDM ABM that have been averaged over different numbers of ABM simulations. All data in this case study is simulated using the parameter values $P_p=0.01,P_d=0.005,\text{ and }P_m=1$. ABM data is comprised of $\mean{C_\text{ABM}(t)}$, which is computed over $N=1,5,10,\text{ or }25$ simulations, and we used 10 separate datasets for each value of $N$ to investigate how stochastic fluctuations affect the final results. For DE learning, we use the Greedy algorithm to solve $\text{d}C/\text{d}t=\bm{\Theta}\xi$ for $\bm{\Theta}=[C,C^2,C^3,C^4]$. We denote $\hat{\xi}^N$ as the estimate of $\xi$ that results for each value of $N$. The code for this case study is provided in the file {\color{blue} \textbf{Case study 2 varying number of ABM Simulations.ipynb}}.

\begin{figure}
    \centering
    \includegraphics[width=0.99\textwidth]{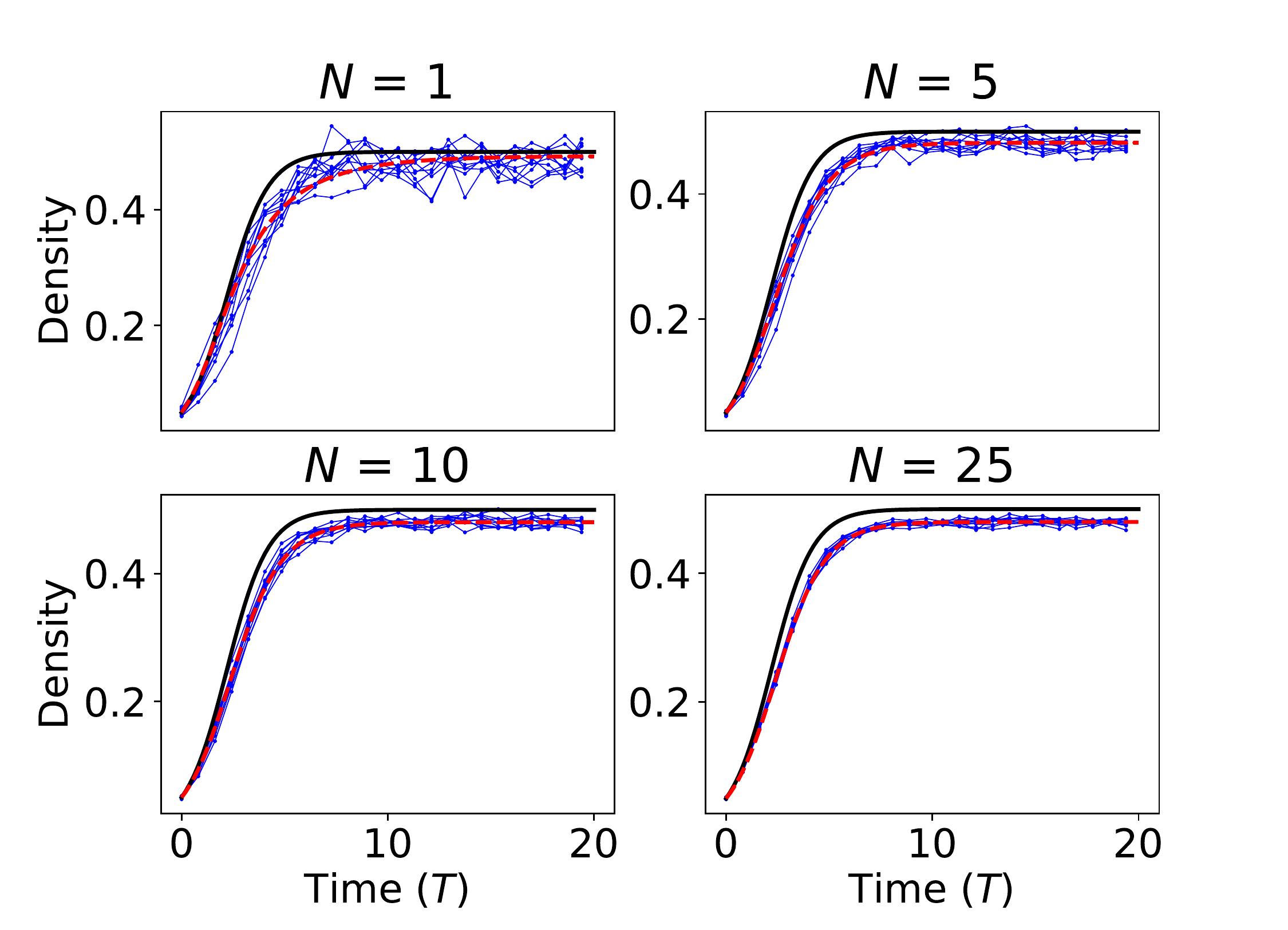}
    \caption{Case study 2. Learning equations from varying numbers of ABM simulations.  In each figure, we depict ten realisations of $\mean{C_\text{ABM}(t)}$ (blue dots and line) against the corresponding mean-field model (solid black line) and averaged learned equation, (red dashed line), which depicts the final learned equation whose coefficients were averaged over all the learned DE models for each realisation of $\mean{C_\text{ABM}(t)}$.
    Each realisation of $\mean{C_\text{ABM}(t)}$ was averaged over $N=1$ (top left), $N=5$ (top right), $N=10$ (bottom left), and $N=25$ (bottom right) simulations of the BDM model. \new{Each individual simulation was initialized by placing agents uniformly at random throughout the lattice so that 5\% of lattice sites were occupied.} All simulations are depicted against nondimensionalised time $T=t(P_p-P_d)$. The ABM parameters used in these simulations are
    $P_p=0.01,P_d=0.005,P_m=1$.}
    \label{fig:CS_2_real_convergence}
\end{figure}

For each value of $N$ considered, we learned ten separate DE models from ten separate realisations of $\mean{C_\text{ABM}(t)}$ and then computed the average learned DE model by averaging the coefficients from each of the ten learned equations. Model predictions from the averaged learned DE model and mean-field model are depicted against all ABM data in Figure \ref{fig:CS_2_real_convergence}. The solution profile for each average learned DE model does not change too much as $N$ increases, but the ABM data fluctuations decrease with larger values of $N$. We depict the distributions of each coefficient for each value of $N$ in Figure \ref{fig:app_bw} in Appendix \ref{sec:bw}. As expected, the variance of each distribution decreases with $N$, so we can be more certain about learned DE models with more ABM simulations. We present the averaged learned DE models for each value of $N$ in Table \ref{tab:CS_2_BDM} as well as the averaged mean-field and learned model MSEs. The mean-field MSE maintains a nearly constant value for all $N$, but the learned model MSE decreases with $N$. The learned equation is cubic for all values of $N$, and the cubic coefficient appears to approach zero as $N$ increases. The differences in MSE between successive averaged $\xi^N$ estimates (i.e., $\| \hat{\xi}^5 - \hat{\xi}^1 \|_2$, etc.) is low for $N>10$. The insensitivity of the learned equation above $N=10$ suggests that $N=10$ or 25 simulations are sufficient to accurately capture the mean BDM ABM dynamics for the parameter values we used. \new{As discussed in Section \ref{subsubsec:EQL_pipeline}, we used finite differences for numerical differentiation in these cases. Recent studies have demonstrated how the use of polynomial spline interpolation or artificial neural networks (ANNs) to improve EQL performance in the presence of noise levels on the magnitude of that observed in experimental data \citep{lagergren_learning_2020,rudy_data-driven_2017}}.

\begin{table}

    \begin{tabular}{|c|l|l|c|}
    \hline
    $N$  & Mean-Field Model (MSE) & Learned Model (MSE) & $\hat{\xi}$ MSE \\ 
    \hline
    1 & $\nicefrac{\text{d}C}{\text{d}t} = 0.005C    - 0.01C^2   $ (0.0037) & $\nicefrac{\text{d}C}{\text{d}t} = 0.00568C    - 0.01953C^2    + 0.01624C^3   $ (0.0025) & - \\
    \hline
    5 & $\nicefrac{\text{d}C}{\text{d}t} = 0.005C    - 0.01C^2   $ (0.0031) & $\nicefrac{\text{d}C}{\text{d}t} = 0.00482C    - 0.01299C^2    + 0.00622C^3   $ (0.0012) & 0.012 \\
    \hline
    10 & $\nicefrac{\text{d}C}{\text{d}t} = 0.005C    - 0.01C^2   $ (0.0028) & $\nicefrac{\text{d}C}{\text{d}t} = 0.00472C    - 0.01193C^2    + 0.00439C^3   $ (0.0008) & 0.002 \\
    \hline
    25 & $\nicefrac{\text{d}C}{\text{d}t} = 0.005C    - 0.01C^2   $ (0.0027) & $\nicefrac{\text{d}C}{\text{d}t} = 0.00453C    - 0.01054C^2    + 0.00232C^3   $ (0.0005) & 0.003 \\
    \hline
    \end{tabular}

    \caption{Case study 2. Learned DE models for the BDM process for various numbers of ABM simulations. We fixed $P_p=0.01,P_d=0.005,\text{ and }P_m=1$ in each scenario and averaged ABM output over the given value of $N$ ABM simulations ten separate times to investigate the EQL method's performance in the presence of stochastic ABM fluctuations. The presented learned DE models depicts the final learned equation whose coefficients were averaged over all the learned DE models for each realisation of $\mean{C_\text{ABM}(t)}$. The rightmost column corresponds to the MSE between successive $\hat{\xi}$ estimates: e.g., for $N=5$, we compute $\| \hat{\xi}^5 - \hat{\xi}^1 \|_2.$}
    \label{tab:CS_2_BDM}
\end{table}

\subsection{Case study 3: Learning ABM dynamics from sparse time samples} \label{subsec:case_study_3_unobserved}

A current challenge for modellers is to develop EQL methods that are able to learn DE models from real experimental or clinical data. ABMs are a useful intermediate step to test the predictions of mathematical methods because ABMs emulate the stochastic and discrete nature of many biological processes and allow researchers to alter aspects of the data. Biological data presents many practical challenges for modellers, including only partial observations of the process under consideration or sparse sampling of the data \citep{nardini_learning_2020}. We will use this case study to consider the performance of the EQL methods in the face of both limited data sampling and partial data observations. In turn, we address the following questions: 
\begin{displayquote}
(Q4) How can we determine the resolution needed for accurate DE model learning?\\
(Q5) Which time scales are informative for learning predictive DE models for unobserved data?
\end{displayquote}

\subsubsection{Case study 3a: Learning BDM dynamics from sparsely sampled ABM data} 

We applied the EQL methodology to ABM data where $n=13,25,50\text{ and }100$ time samples were collected. For all values of $n$, the data were collected at equispaced time intervals between the same starting and ending time points. All ABM simulations were computed with mechanistic parameters $P_p=0.05,P_d=0.0125,$ and $P_m=1$. The ABM data, $\mean{C_\text{ABM}(t)}$, was averaged over $N=50$ simulations. The code for this case study is provided in the file {\color{blue} \textbf{Case study 3a Varying data resolution.ipynb}}.

\begin{figure}
    \centering
    \includegraphics[width=.99\textwidth]{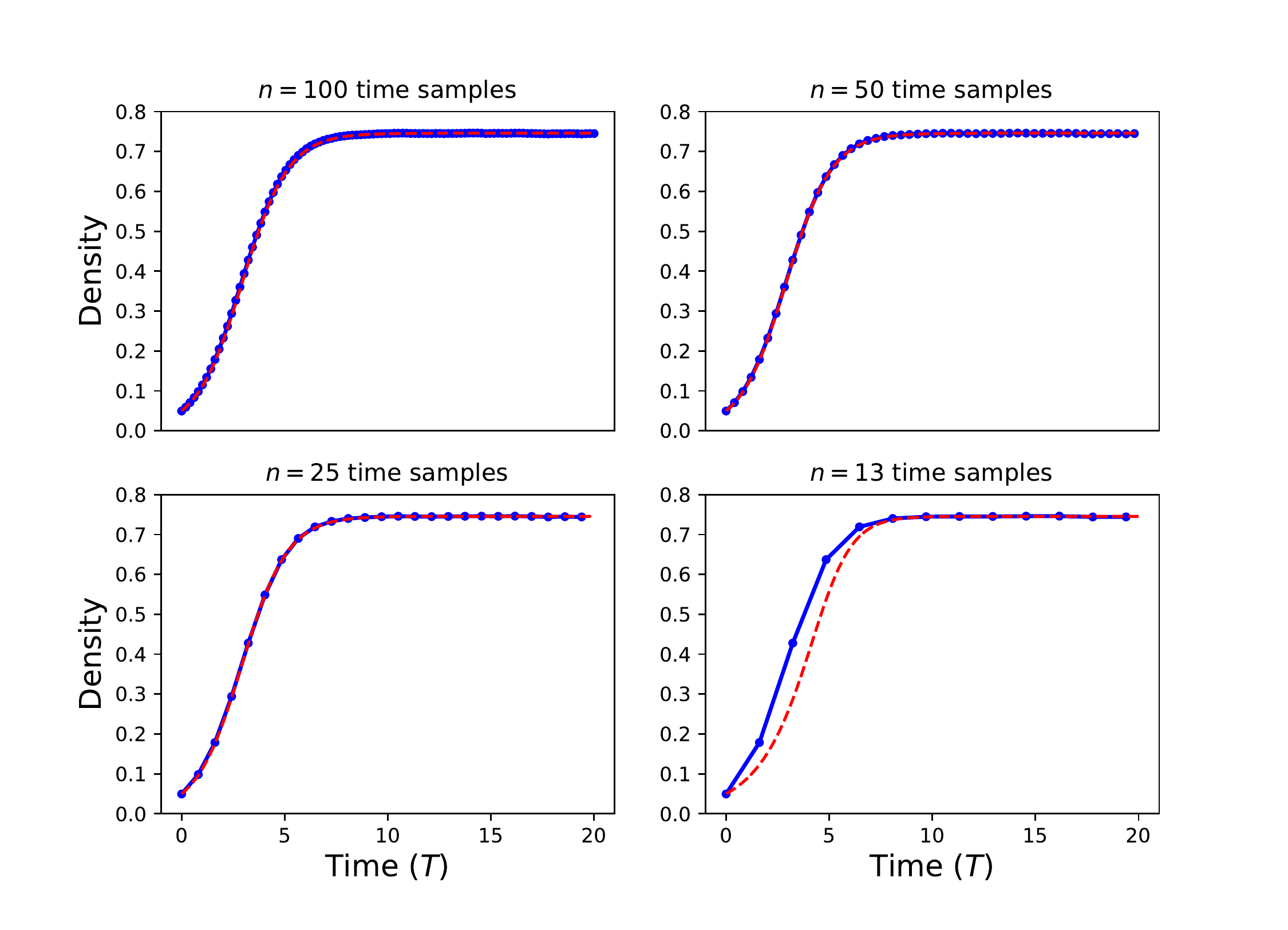}
    \caption{Case study 3a. Learning equations with varying time resolution. We applied the EQL methodology for ABM data with $n=13,25,50,$ or 100 time samples and depict the learned model (red dashed line) against the ABM data (blue solid line and dots). All simulations are depicted against nondimensionalised time $T=t(P_p-P_d)$. We fixed $P_p=0.05,P_d=0.0125,\text{ and }P_m=1$.}
    \label{fig:CS_3_BDM}
\end{figure}

We depict model predictions from the learned equation against ABM data in Figure \ref{fig:CS_3_BDM}. The learned DE models accurately predict the ABM output for $n=100, 50, \text{ and } 25$ time samples. With $n=13$ time samples, however, the learned equation predicts the same carrying capacity as the data but fails to accurately predict the ABM dynamics before plateauing (excluding the initial time sample, which is used as the initial condition). The learned model equations for these different scenarios, and the MSE between the learned model prediction and ABM data from each model, are summarized in Table \ref{tab:CS_3_BDM}. The MSE is less than or equal to $10^{-3}$ for $n\ge13$ but rises to 0.0154 for $n=0.13$. From this case study, we suggest that $n=25$ time samples provides sufficient time resolution to accurately infer the underlying dynamics for parameter values $P_p=0.05$, $P_d=0.0125$, and $P_m=1$. If we analysed experimental data we had reason to believe resulted from similar underlying parameter values, we would trust our analysis of the experimental data if we had 25 or more equispaced time points over this same time interval. \new{We note alternative methods of numerical differentiation can be used to improve EQL results with limited time samples here. As an example, we have recently introduced an ANN approach that has proven successful in learning equations from sparse PDE data \citep{nardini_learning_2020}. This method would likely outperform the finite difference scheme used for differentiation throughout this study. We point the interested reader to \citep{lagergren_learning_2020,nardini_learning_2020} for more information on this ANN approach.}

\begin{table}[]
    \centering
    \begin{tabular}{|c|l|}
    \hline
    \% of data & \multicolumn{1}{c|}{Learned Model (MSE)} \\ 
    \hline
    $100$ & $\nicefrac{\text{d}C}{\text{d}t} = 0.03374C    - 0.04522C^2   $ (0.0002) \\
    \hline
    $50$ & $\nicefrac{\text{d}C}{\text{d}t} = 0.03379C    - 0.04531C^2    - 0.0C^3   $ (0.0002) \\
    \hline
    $25$ & $\nicefrac{\text{d}C}{\text{d}t} = 0.03372C    - 0.0452C^2   $ (0.0004) \\
    \hline
    $13$ & $\nicefrac{\text{d}C}{\text{d}t} = 0.02073C    - 0.03733C^3   $ (0.0154) \\
    \hline
    \end{tabular}
    \caption{Case study 3a. Learned model equations for case study 3a with varying numbers of time samples, $n$, from the data. MSE denotes the mean squared error between the learned model prediction and $\mean{C_\text{ABM}(t)}$. }
    \label{tab:CS_3_BDM}
\end{table}

\subsubsection{Case study 3b: Predicting unobserved BDM dynamics} 

We applied the EQL methodology to ABM data where only the first 10\%, 20\%, 25\% or 50\% of the ABM data are used for training and the remaining data are held out for testing how well the learned DE model predicts unobserved dynamics. ABM data is computed with $P_p=0.01,P_d=0.005, \text{ and }P_m=1$ and $\mean{C_\text{ABM}(t)}$ is averaged over 50 ABM simulations with $n=100$ data points until $t=15(P_p-P_d)$. The code for this case study is provided in the file {\color{blue} \textbf{Case study 3b predicting unobserved dynamics.ipynb}}.

\begin{figure}
    \centering
    \includegraphics[width=.99\textwidth]{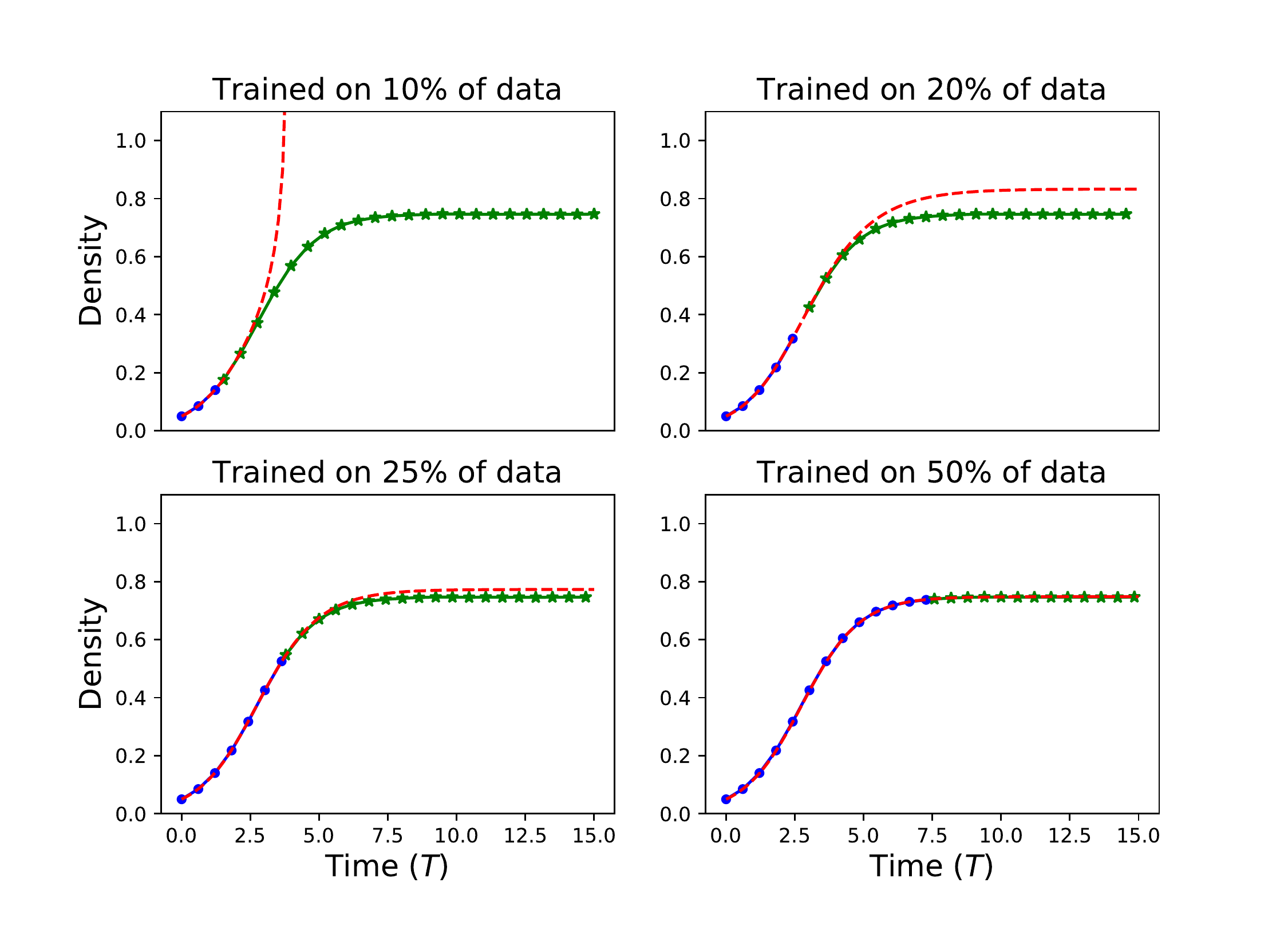}
    \caption{Case study 3b. Predicting BDM dynamics from partial ABM data. We applied the EQL methodology on the first 10\% (top left), 20\% (top right), 25\% (bottom left), and 50\% (bottom right) of $\mean{C_\text{ABM}(t)}$ to investigate the learned equations performance in predicting unoberved ABM dynamics. The blue dots correspond to ABM output that was used to infer the learned model, the green stars denote ABM output that was used for model testing, and the red dashed line denotes the solution of the learned model. All simulations are depicted against nondimensionalised time $T=t(P_p-P_d)$. We fixed $P_p=0.01,P_d=0.005,\text{ and }P_m=1$ in each scenario and averaged ABM output over $N=50$ simulations.}
    \label{fig:CS_3b_BDM}
\end{figure}

We depict model predictions against the testing and training data in Figure \ref{fig:CS_3b_BDM}. When trained on 10\% of data, the inferred model grows without bound for all simulated time and does not accurately describe the testing data. When trained on 20\% of data, the learned model plateaus above the carrying capacity of the data, approximately $C=0.75$. For training data comprised of 25\% and 50\% of the available data, the learned models closely predict the test data. We suggest that, for this case study, data should be sampled beyond the inflection point in order to accurately predict the unobserved dynamics and carrying capacity. \new{One possible explanation for this observation is that all data before the inflection point is concave up and all data after the inflection point is concave down. Including data past the inflection point appears necessary to capture how $C_\text{ABM}(t)$ becomes concave down before reaching the carrying capacity. Further investigation into informative datasets for EQL training have been explored elsewhere \citep{nardini_learning_2020}.} 

The learned DE model for these different data percentages and their MSEs on the test data set are summarized in Table \ref{tab:CS_3b_BDM}. As expected, the testing data MSE decreases as more of the data is used to learn the DE model. 

\begin{table}[]
    \centering
    \begin{tabular}{|c|l|}
    \hline
    \% & \multicolumn{1}{c|}{Learned Model (Testing MSE)} \\ 
    \hline
    $10$ & $\nicefrac{\text{d}C}{\text{d}t} = 0.00802C    - 0.02098C^2    + 0.03837C^3   $ (nan) \\
    \hline
    $20$ & $\nicefrac{\text{d}C}{\text{d}t} = 0.00745C    - 0.01154C^2    + 0.00305C^3   $ (0.0071) \\
    \hline
    $25$ & $\nicefrac{\text{d}C}{\text{d}t} = 0.00737C    - 0.01081C^2    + 0.00159C^3   $ (0.0021) \\
    \hline
    $50$ & $\nicefrac{\text{d}C}{\text{d}t} = 0.00721C    - 0.00975C^2    + 0.00016C^3   $ (0.0002) \\
    \hline
    \end{tabular}
    \caption{Case study 3b. Summary of models learned with the first 10\%, 20\%, 25\%, and 50\% of $\mean{C_\text{ABM}(t)}$.  MSE on the remaining data is given in parentheses. Data were simulated with the BDM model using $P_p=0.01,P_d=0.005,\text{ and }P_m=1$, averaged over $N=50$ ABM simulations. Note that our implementation of the learned DE model trained on 10\% of the data fails to converge because this model grows faster than exponential growth, which is not realistic of the ABM data. We depict the MSE of this learned equation against $\mean{C_\text{ABM}(t)}$ as nan (not a number) due to this numerical instability.}
    \label{tab:CS_3b_BDM}
\end{table}

\subsection{Case study 4: Using EQL methods for model selection} \label{subsec:case_study_4_model_selection}

Our first case study used EQL methods to demonstrate that the mean-field model may not be valid for predicting the dynamics of the BDM ABM when $P_p>0.1$ (with other parameters fixed at $P_m=1\text{ and }P_d = P_p/2$). If one wants to use a DE model to analyse such ABM simulations, then an alternative model may be required. In Appendix \ref{sec:app_deriving_DE}, we show the DE model given by
\begin{equation}
    \dfrac{\text{d}C}{\text{d}t} = P_pC(1-FC) - P_dC,
    \label{eq:modified_reminder}
\end{equation}
can be used to model output from the BDM ABM. In Equation \eqref{eq:modified_reminder}, $F$ is the occupancy correlation between neighbouring lattice sites \citep{baker_correcting_2010}. For the lattice-based BDM model, this value is defined between two neighbouring lattice sites $\alpha,\beta$ as
\begin{equation}
    F(t) = \dfrac{\mathbb{P}[A_\alpha(t),A_\beta(t)]}{\mathbb{P}[A_\alpha(t)]\mathbb{P}[A_\beta(t)]}. \label{eq:corr_occupancy}
\end{equation}
Note that if the occupancy probabilities of these sites are independent, then $F(t)\equiv 1$, and Equation \eqref{eq:modified_reminder} simplifies to the mean-field DE model. 

Model selection studies \citep{bortz_model_2006} are suited to determine which model most parsimoniously describes a given dataset from several plausible models. We may now be interested in determining which of our our two models, the mean-field model in Equation \eqref{eq:logistic}, or the DE model in Equation \eqref{eq:modified_reminder}, best describes ABM output for the BDM process.  A typical model selection study for these two models may be problematic, however, as deriving and computing the DE model for $F$ in Equation \eqref{eq:modified_reminder} is complicated even for a scenario as simple as the BDM model, yielding an additional set of auxiliary DEs needed to describe $F$ \citep{baker_correcting_2010}. Since the derivation of such auxiliary equations is effectively intractable for more complex ABM models, we will investigate the following question:

\begin{displayquote}
(Q6) Can methods from EQL be used for DE model selection for ABM analysis?
\end{displayquote}
In doing so, we propose an alternative strategy for model selection, i.e., for selecting additional time-dependent variables, such as the occupation correlation, $F(t)$, to increase DE model accuracy with concepts from EQL.

\subsubsection{Model selection with EQL for the BDM process}

 We use this section to demonstrate how EQL methods can be used for model selection between the logistic equation given by Equation \eqref{eq:logistic} and the modified logistic equation given by Equation \eqref{eq:modified_reminder}. The code for this case study is provided in the file {\color{blue} \textbf{Case study 4 model selection with EQL.ipynb}}. The rate of proliferation, $P_p$, varies over the values $P_p = 0.005,0.01,0.05,0.1,0.5$. For each value of $P_p$, we set $P_d=P_p/2$ and fix $P_m=1$. The occupancy correlation value from Equation \eqref{eq:corr_occupancy} is estimated from the $n^\text{th}$ of $N$ ABM simulations by computing
\begin{equation}
    F^{(n)}(t) = \dfrac{C_{ABM}^{(n,2)}(t) X^4}{\chi^{2}(C_{ABM}^{(n)}(t))^2},
\end{equation}
where $C_{ABM}^{(n,2)}(t)$ is the number of jointly occupied neighbouring pairs of lattice sites over time, $X$ is the length of the lattice, and $\chi^{2}$ is the total number of adjacent lattice-site pairs \citep{fadai_accurate_2019}. The average occupancy correlation  values is then averaged over all $N=50$ simulations:
\begin{equation}
    \mean{F(t)} = \dfrac{1}{N}\sum_{k=1}^N F^{(n)}(1,t).
\end{equation}
 
 The model selection methodology using EQL concepts proceeds as follows. From each ABM dataset, we compute both $C_d=\mean{C_\text{ABM}(t)}$ and $F_d=\mean{F(t)}$ and substitute these values into two separate $n\times 2$ matrices of right hand side terms given by
\[
\bm{\Theta}_{1}=[C_d(1-C_d),C_d],\ \ \ \bm{\Theta}_{2}=[C_d(1-FC_d),C_d],
\]
where any multiplication is performed element-wise. Note that $\bm{\Theta}_{1}$ corresponds to the library of terms for Equation \eqref{eq:logistic} while $\bm{\Theta}_{2}$ corresponds to the library of terms for Equation \eqref{eq:modified_reminder}. We uniformly at random place half of the elements comprising $\text{d}C_d(t)/\text{d}t$ and the corresponding rows from $\bm{\Theta}_{1},\bm{\Theta}_{2}$ into training sets given by the vector $\text{d}C^\text{train}_d(t)/\text{d}t$ and the $n/2\times 2$ matrices $\bm{\Theta}^\text{train}_{1},\bm{\Theta}^\text{train}_{2}$. The remaining elements of $\text{d}C_d(t)/\text{d}t$ and rows of $\bm{\Theta}_{1},\bm{\Theta}_{2}$ are placed into testing sets given by $\text{d}C^\text{test}_d(t)/\text{d}t$ and matrices $\bm{\Theta}^\text{test}_{1},\bm{\Theta}^\text{test}_{2}$. The training set is used to estimate $\hat{\xi}_1$ from $\bm{\Theta}^\text{train}_{1}$ and $\hat{\xi}_2$ from $\bm{\Theta}^\text{train}_{2}$
by solving the two linear regression problems 
\begin{equation}
\dfrac{\text{d} C^\text{train}_d}{\text{d}t}=\bm{\Theta}^\text{train}_{i}\xi_{i},\ i=1,2.\label{eq:ODE_find}
\end{equation}
Note that the vector $\hat{\xi}_1$ parameterizes Equation \eqref{eq:logistic}, and $\hat{\xi}_2$ parameterizes Equation \eqref{eq:modified_reminder}. We then use the
testing set to select the best model for each data set by \begin{equation}
 \hat{\xi}_i = \arg\min_{i}\left\|\dfrac{\text{d}C^\text{test}_d}{\text{d}t}-\bm{\Theta}^\text{test}_{i}\hat{\xi}_{i}\right\|_{2}.\label{eq:EQL_model_selection}
\end{equation} 
We use 100 randomly sampled training and validation sets and select whichever of the two models minimizes Equation \eqref{eq:EQL_model_selection} more often in these 100 testing-validation realisations.

In Table \ref{tab:CS_4_model_selection}, we present the final selected models for various values of $P_p$. The mean-field model is selected for $P_p=0.005 \text{ and }0.01$ with 57 and 77 of the 100 total votes, respectively. Equation \eqref{eq:modified_reminder} is selected for all larger proliferation values with at least 93 of the 100 total votes. These results are in agreement with Case study 1, where the mean-field model predicted ABM output well at $P_p =0.01$, but was unable to predict these data for larger values of $P_p$ (see, e.g., Figure \ref{fig:CS_1_BMD_data_comparison}). In addition to more accurately matching the ABM output than the mean-field model, Equation \eqref{eq:modified_reminder} also provides accurate parameter estimates for $P_p$ and $P_d$. By matching the forms of the modified logistic equation to the learned equation (i.e., we can estimate $P_p$ using the coefficient in front of $C(1-FC)$ and $P_d$ using the negative of the coefficient in front of $C$), we observe that these estimates appear very close to their true underlying values.

\begin{table}
\centering
\begin{tabular}{|c|c|c|c|}
    \hline
    $P_p$ & \ $P_d$ & \ Selected Model & Votes (out of 100) \\ 
    \hline
    $0.005$ & $0.0025$ & $\nicefrac{\text{d}C}{\text{d}t} =  - 0.00245C    + 0.00485C(1-C)   $ & 57 \\
    \hline
    $0.01$ & $0.005$ & $\nicefrac{\text{d}C}{\text{d}t} =  - 0.00483C    + 0.00952C(1-C)   $ & 77 \\
    \hline
    $0.05$ & $0.025$ & $\nicefrac{\text{d}C}{\text{d}t} =  - 0.02482C    + 0.04936C(1-FC)   $ & 93 \\
    \hline
    $0.1$ & $0.05$ & $\nicefrac{\text{d}C}{\text{d}t} =  - 0.04966C    + 0.09874C(1-FC)   $ & 100 \\
    \hline
    $0.5$ & $0.25$ & $\nicefrac{\text{d}C}{\text{d}t} =  - 0.25248C    + 0.50271C(1-FC)   $ & 100 \\
    \hline
\end{tabular}
\caption{Case study 4. Model selection with EQL. We present the average selected model equations for $\langle C_\text{ABM}(t)\rangle$ over various values of $P_{p}$. For each value of $P_p$, we set $P_{d}=P_p/2$ and $P_{m}=1$. ABM output is averaged over $N=50$ ABM simulations of the BDM model to ensure convergence to mean behaviour. The rightmost column lists how many votes the selected equation received out of 100 total.\label{tab:CS_4_model_selection}}
\end{table}

\subsection{Case study 5: Varying infection probability rates for the SIR model}\label{subsec:varying_infection_SIR}

 The EQL methods considered in this work are applicable to many ABMs, including the SIR model introduced in Section \ref{sec:SIR_ABM}. We will now learn models for this ABM over several parameter regimes and test the performance of the mean-field and learned models in predicting ABM output \citep{blackwood_introduction_2018,diekmann_definition_1990, viceconte_covid-19_2020}.

We let the agent infection rate take the values $P_I=$ $0.005,0.01,0.05,0.1$, set the agent recovery rate, $P_R$, to be one-tenth of the infection rate for each value of $P_I$, and fix  $P_m=1$ for all scenarios. The ABM output is comprised of $S_d = \mean{S_\text{ABM}(t)}$, $I_d = \mean{I_\text{ABM}(t)}$, and $R_d = \mean{R_\text{ABM}(t)}$, all of which are averaged over $N=25$ ABM simulations from a square lattice with length $X=40$. Because these three quantities will always sum to unity in the model, we focus on learning equations for $S(t)$ and $I(t)$ and note that $R(t)=1-S(t)-I(t)$. For model learning, we use the matrix of potential right-hand side terms given by $\bm{\Theta} = [S,S^2, I,I^2,SI]$. The Lasso algorithm is used here to solve for the unknown vectors $\xi_1$ and $\xi_2$ from the linear systems $\text{d}S/\text{d}t = \bm{\Theta}\xi_1$ and $\text{d}I/\text{d}t = \bm{\Theta}\xi_2$. Note that even for a two-dimensional system with five library terms, the total number of possible models for the $S$ and $I$ variables is $\sum_{i=0}^5 {5 \choose i}=32$ each, resulting in a total combination of $32^2=1024$ possible DE models, which highlights the difficulty in finding a predictive model from this large suite of possibilities. The code for this case study is provided in the file  {\color{blue} \textbf{Case study 5 SIR Varying params.ipynb}}.

\begin{figure}
    \centering
    \includegraphics[width=0.9\textwidth]{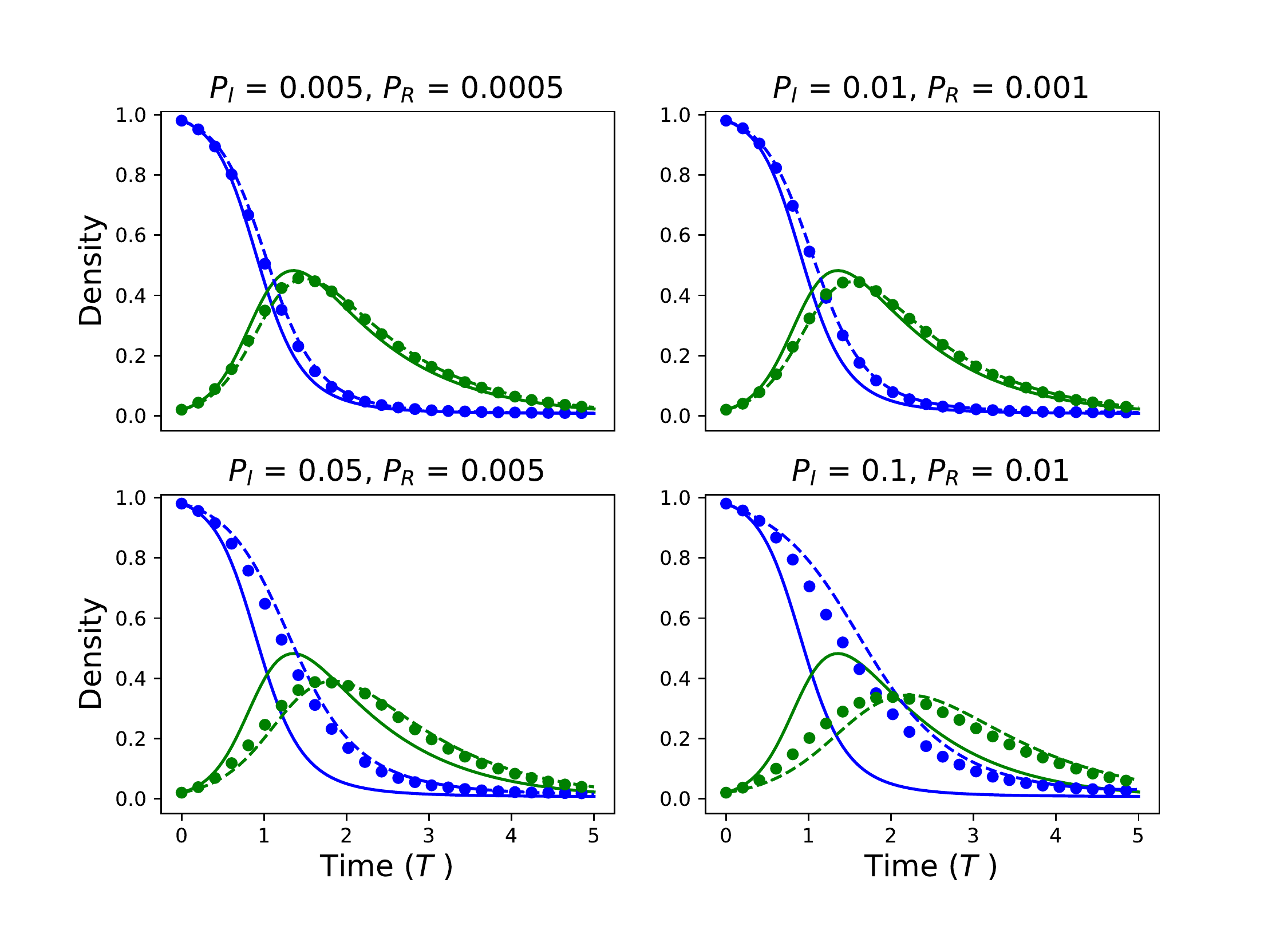}
    \caption{Case study 5. Comparing mean-field and learned model predictions to ABM data from the SIR model. In each figure, we depict the predictions of the mean-field model (solid blue curve for $S(t)$ and solid green curve for $I(t)$), against the predictions of the learned DE models (blue dashed curve for $S(t)$ and dashed green line for $I(t)$), as well as $\mean{S_\text{ABM}(t)}$ (blue dots) and $\mean{I_\text{ABM}(T)}$ (green dots). All simulations are depicted against nondimensionalised time $T=tP_R$.  We fixed $P_m=1$ for each simulation. All ABM simulations were computed on a square lattice of length $X=40$.}
    \label{fig:case_study_1_SIR_compare}
\end{figure}

We depict the predictions of the mean-field and learned DE models over time against ABM output in Figure
\ref{fig:case_study_1_SIR_compare}. The mean-field and learned DE model equations with corresponding MSEs are presented in Table \ref{tab:SIR_learned_eqs}. For $P_I=0.005$ the mean-field model predicts the ABM output well, but as $P_I$ increases the mean-field model predictions worsen, as evidenced by increases in the MSE. The mean-field model underpredicts the susceptible agent density at all times and overpredicts the infected agent density at early times.  The learned model, on the other hand, is able to predict ABM output accurately for all values considered and achieves low MSE values (with the lowest value at $P_I=0.01$). The learned equation forms are similar to the mean-field model for $P_I=0.005$, 0.01, and 0.05. For $P_I=0.1,$ the learned equations have additional $S,S^2,$ and $I$ terms in the model equations for $\text{d}S/\text{d}t$. These terms may be used to capture effects that are not described by the mean-field model. 

The basic reproductive number, or $R_0$,  is defined as the expected number of secondary infections that result from a single infection in a population comprised solely of susceptible agents \citep{blackwood_introduction_2018,diekmann_definition_1990,viceconte_covid-19_2020}. A disease may spread rapidly when $R_0>1$, and will die out when $R_0<1$. $R_0$ is now a commonly-used criterion to determine when a disease will continue to spread through a population and possibly cause an outbreak. More details of how to computer $R_0$ are provided in \citep{blackwood_introduction_2018}, but we note here that $R_0$ can be found for the mean-field SIR model by determining critical parameter values where d$I(t=0)/$d$t$ will exceed zero at the initial condition. From Equation \eqref{eq:SIR_MF}, the mean-field model shows that d$I/$d$t>0$ when $MP_IS>P_R$. Recalling that $S(t)\approx1$ at the start of the disease spread, this implies that d$I/$d$t>$0 when $MP_I/P_R>1$, i.e., $R_0 = MP_I/P_R$ for the mean-field model. The values of $R_0$ can be computed using similar methods for the learned equations detailed in Table \ref{tab:SIR_learned_eqs}. Here we observe that the mean-field model predicts that $R_0=5.0$ for all considered scenarios (recall, $M=0.5$ in all simulations), yet as $P_I$ increases from 0.005 to 0.01,0.1,0.25 (with $P_I/P_R=10.0$), the learned equations predict $R_0 = 4.68,4.52,3.72, 3.41$, respectively. Thus, while the mean-field model suggests that a single infected individual will cause five secondary infections (in a population full of susceptible agents) at each of these infection rates, the learned models suggest that as $P_I$ increases (with the ratios $P_I/P_R$ and $P_I/P_m$ fixed), the number of secondary infections will decrease in this scenario. Recall from Figure \ref{fig:SIR_ABM_compare} that infected agents tend to cluster with large values of $P_I$ (while $P_m$ and $P_I/P_R$ are fixed). The mean-field model does not consider how such spatial clustering may affect the number of secondary infections that result from a single infection in a population comprised solely of susceptible agent. The learned DE models, on the other hand, implicitly account for this information when finding a predictive DE model from the ABM data.

\begin{landscape}
\begin{table}
    \begin{tabular}{|c|c|l|c|l|c|}
    \hline
    $P_I$ & \ $P_R$ & \ Mean-field model (MSE) & $R_0$ & Learned model (MSE) & $R_0$ \\ 
    \hline
    \multirow{2}{*}{$0.005$} & \multirow{2}{*}{$0.0005$}  & $\nicefrac{\text{d}S}{\text{d}t} =  - 0.0025IS   $ (0.0027) & \multirow{2}{*}{5.0}  & $\nicefrac{\text{d}S}{\text{d}t} =     - 0.00229IS   $ (0.0012) & \multirow{2}{*}{4.68} \\
     & & $\nicefrac{\text{d}I}{\text{d}t} = 0.0025IS    - 0.0005I   $ (0.002) & & $\nicefrac{\text{d}I}{\text{d}t} =  - 0.00049I    + 0.00225IS   $ (0.0009) & \\
    \hline
    \multirow{2}{*}{$0.01$} & \multirow{2}{*}{$0.001$}  & $\nicefrac{\text{d}S}{\text{d}t} =  - 0.005IS   $ (0.0044) & \multirow{2}{*}{5.0}  & $\nicefrac{\text{d}S}{\text{d}t} =   - 0.00447IS   $ (0.0006) & \multirow{2}{*}{4.52} \\
     & & $\nicefrac{\text{d}I}{\text{d}t} = 0.005IS    - 0.001I   $ (0.0029) & & $\nicefrac{\text{d}I}{\text{d}t} =  - 0.00098I    + 0.00443IS   $ (0.0005) & \\
    \hline
    \multirow{2}{*}{$0.05$} & \multirow{2}{*}{$0.005$}  & $\nicefrac{\text{d}S}{\text{d}t} =  - 0.025IS   $ (0.0107) & \multirow{2}{*}{5.0}  & $\nicefrac{\text{d}S}{\text{d}t} =  - 0.01881IS   $ (0.003) & \multirow{2}{*}{3.72} \\
     & & $\nicefrac{\text{d}I}{\text{d}t} = 0.025IS    - 0.005I   $ (0.0066) & & $\nicefrac{\text{d}I}{\text{d}t} =  - 0.00472I    + 0.01833IS   $ (0.0021) & \\
    \hline
    \multirow{2}{*}{$0.1$} & \multirow{2}{*}{$0.01$}  & $\nicefrac{\text{d}S}{\text{d}t} =  - 0.05IS   $ (0.0164) & \multirow{2}{*}{5.0}  & $\nicefrac{\text{d}S}{\text{d}t} = 0.01569S    - 0.01608S^2    - 0.00545I    - 0.0459IS   $ (0.0055) & \multirow{2}{*}{3.41} \\
     & & $\nicefrac{\text{d}I}{\text{d}t} = 0.05IS    - 0.01I   $ (0.0097) & & $\nicefrac{\text{d}I}{\text{d}t} =  - 0.00906I    + 0.03101IS   $ (0.0033) & \\
    \hline
\end{tabular}
    \caption{Case study 5: Mean-field and learned DE models for the SIR ABM for various values of $(P_I,P_R)$, along with the computed MSE  values between model and ABM output, and model $R_0$ calculations. We fix $P_m=1$ in each scenario.}
    \label{tab:SIR_learned_eqs}
\end{table}
\end{landscape}

\section{Conclusions and discussion} \label{sec:conclusions}

In this work, we have considered three different, yet synergistic, approaches to study the emergent behaviour of ABMs. These approaches include extensive simulation, DE model derivation using coarse graining approaches, and EQL. We demonstrated some of the strengths and weaknesses of each approach through their applications to two ABMs: a BDM model and an SIR model. We summarize some of these strengths and weaknesses in Table \ref{tab:pro_cons}.

Extensive simulation is the most commonly-used and straightforward approach in understanding the behaviour of ABMs. Extensive simulation is advantageous because, in theory, it can be used to analyse any ABM: it only requires simulating the ABM on a computer. This approach quickly becomes practically challenging, however, because it can become difficult to simulate sufficiently complex ABMs with basic computing hardware. The BDM process used in this work is simple enough to be simulated on a personal laptop (2.5 GHz Intel Core i7 processor, 16 GB RAM) for all datasets generated for this study. The SIR model is more involved, however, and was run on a computer cluster to ease implementation. Computation of a single SIR simulation with $P_I=0.005,P_R=0.0005,P_M=1.0$ took about 48 minutes to compute. Two further drawbacks of performing extensive simulation are that this technique is not compatible with analytical methods, and that its output may not be interpretable, i.e., we observe that the system behaviour changes as parameters are varied, but we may not understand \emph{why} this happens, nor can we necessarily predict this behaviour \emph{a priori}.

There is a wide literature demonstrating that coarse graining approaches can be used to derive DE models that approximate ABM dynamics \citep{baker_correcting_2010,bernoff_agent-based_2020,binny_living_2020,fadai_accurate_2019,johnston_co-operation_2017,johnston_how_2014,nardini_modeling_2016}. Such DE models are advantageous because they are usually simple to solve (either analytically or numerically), interpretable, and provide insight into how changing ABM parameters can lead to emergent behaviours. Analytical techniques allow us to infer how the system will behave over many different parameter values without simulating the ABM. For example, the per-capita growth rate for the BDM model predicts the population growth rate as a function of density, and the $R_0$ calculation for the SIR model predicts if a disease will outbreak or die out. However, we saw throughout the presented case studies that these analyses fail to accurately predict ABM behaviour for parameter regimes where the mean-field assumption is violated.  Another limitation of this approach is that the derivation of DE models \new{requires user-made assumptions, such as the mean-field assumption, which are often violated by many ABM simulations. As there is no universal methodology to convert ABM rules into predictive DE models for all input parameter values, EQL provides a powerful framework to learn informative DE models for ABMs and, in turn, determine when simplifying  assumptions (\emph{e.g.}, the mean-field assumption) are reasonable.} \old{is only possible for simple ABMs, since it can be challenging to convert complex ABM rules into DE models. }

EQL is a recent field of research that seeks to infer DE models directly from observed data. EQL combines many benefits of the two previously-mentioned methodologies to analyse the emergent behaviour of ABMs. We highlighted many of the advantages of these approaches, and addressed several ways in which EQL can aid modellers in ABM analysis, through our investigation of Questions (Q1)-(Q6). While the mean-field model does not accurately predict ABM output for many parameter combinations, our exploration of (Q1) demonstrated that EQL provides a simple way to determine when the mean-field model can or cannot be trusted for such predictions. If the learned equation form matches the mean-field model, then the mean-field model should provide accurate insight; when it does not, then alternative models may be needed, as we discussed in (Q2). We also demonstrated that EQL methods can be used to infer novel DE models for ABMs. We observed that the mean-field model cannot predict output of the BDM process for large rates of agent proliferation and death relative to motility. Instead, Equation \eqref{eq:learned_CS1}
can be used to accurately predict ABM dynamics. Further investigation needs to be carried out in order to understand how Equation \eqref{eq:learned_CS1} may result from ABM rules.

A significant challenge of using ABMs in practice is their intensive computational nature. Through (Q3), we explored how many ABM simulations are necessary to reliably predict ABM output. We showed that the learned equation may deviate from the average ABM behaviour when only a small number of ABM simulations are used. With a sufficiently large number of ABM simulations, however, the learned equation can accurately predict ABM dynamics. For the BDM process, we found that $N=10$ ABM simulations was sufficient. In practice, we can determine when enough simulations have been performed to capture mean ABM dynamics by considering how much the learned vector, $\xi$, changes with additional ABM simulations. When this vector becomes sufficiently insensitive to increases in $N$, then sufficient ABM dynamics have been performed.

In (Q4), we investigated the sampling resolution needed in time to reliably predict ABM dynamics. In this case study, we determined a magnitude between uniformly-sampled time samples above which the EQL pipeline could not learn a DE model that accurately predicted ABM behaviour. We used Question (Q5) to determine when EQL methods can be used to predict unobserved ABM dynamics. In the presented case study, we used the BDM model to demonstrate that the learned equations can accurately predict unobserved ABM dynamics when the observed data exceeds half of the population's carrying capacity. These two investigations, Q4 and Q5, suggest important future work must be performed to determine strategic (and not necessarily uniform) samples of the ABM that are informative and capture all dynamic regimes of the data). Some preliminary work towards these questions for PDE models has been investigated in \citep{nardini_learning_2020}. A limitation of EQL methods, as opposed to ABM simulation and mean-field models, is that it may not be able to accurately extrapolate to unobserved data and parameters, as we observed in Case study 3b. This is possible with ABM simulation, where the ABM can be run for longer time or at different parameter values. Similarly, the mean-field model can extrapolate its predictions by solving the model over a longer time period or simulating it at different parameter values.

Finally, we considered Question (Q6) to determine if EQL methods can be used to aid in model selection for ABMs. This is advantageous because DE models that are more complex than mean-field models can also be derived for ABMs. Although potentially more accurate, these models are difficult to interpret and simulate than mean-field models. Equation \eqref{eq:modified_reminder}, as one example, can be challenging to implement numerically because the dynamical system for the occupancy correlation function, $F(t)$, requires solving a high-dimensional system of equations using user-defined closure approximations \citep{baker_correcting_2010}. Instead, we can measure $F$ over each ABM simulation and use these observations of $F$ in a hybrid approach to select whether including $F$ leads to a more accurate DE model. Applying such a hybrid approach to the simulation of DE models has been recently proposed to aid in reducing identifiability-related issues \citep{hamilton_hybrid_2017,lagergren_forecasting_2018}. 

EQL methods are quickly growing in popularity as a means to infer DE models from noisy data \citep{lagergren_learning_2020,rudy_data-driven_2017,nardini_learning_2020}. We have shown in this study that such methods provide a reliable and promising tool to aid modellers in interpreting and analyzing ABMs. Learning DE models from data does not require user-made assumptions on whether or not the ABM simulation satisfies certain properties, as is needed for the derivation of mean-field models. \new{We have also demonstrated in this work how EQL methods can be used to} predict ABM dynamics from limited ABM simulations, \new{learn DE models from only a subset of data, and }accurately predict dynamics over a wide range of parameter values. \new{Such learned equations from ABM data also make ABM analysis more interpretable, as analysis of the learned equation provides insight into the underlying biological mechanisms (\emph{e.g.}, per-capita growth rates, $R_0$, bifurcation analysis, etc.). The order of a learned equation may provide insight into how many neighboring lattice site occupancies impact individual agent behaviour \citep{simpson_cell_2010}. }

\new{There are many areas for exciting future work in ABM analysis using EQL methods. For example, global sensitivity analysis techniques \citep{smith_uncertainty_2013} could be used to determine ABM parameter thresholds where the learned equation forms change, and what insights these threshold values may provide into the ABM dynamics.} We anticipate that, through this tutorial, EQL will increasingly be used to interpret complex ABM simulations. \new{Future work should aim to address challenges that prevent the learning of DE models from real experimental data. While we focused on learning deterministic DE models from stochastic ABM simulations in this work, more recent studies have explored learning stochastic DE model forms (including both drift and diffusion estimates) from data \citep{klus_data-driven_2020}. Another recent study has shown that the dynamics from a stochastic non-Markovian model can be learned using a simpler time-inhomogenous Markovian model framework with the aid of ANNs \citep{jiang_neural_2020}. We also focused on simple ABMs in this study (including the BDM process and the SIR model), but future work should examine model learning for more complex ABM dynamics, such as bistable \citep{johnston_co-operation_2017} and periodic behaviour \citep{walker_agent-based_2004}. }

\begin{table}[]
    \centering
    \begin{tabular}{|c|c|c|c|}
    \hline
     & Extensive simulation & DE Model derivation & Equation learning \\     
    \hline
    Amenable for ABMs  & \multirow{2}{*}{\checkmarkJTN} & \multirow{2}{*}{\xmark} & \multirow{2}{*}{\checkmarkJTN}  \\
    of varying complexity & & & \\
    \hline
    Computationally efficient & \xmark & \checkmarkJTN & \checkmarkJTN  \\
    \hline
    Amenable to analytical  & \multirow{2}{*}{\xmark} & \multirow{2}{*}{\checkmarkJTN} & \multirow{2}{*}{\checkmarkJTN}  \\
    forms of exploration & & & \\
    \hline
    Provides accurate estimates  & \multirow{3}{*}{\checkmarkJTN} & \multirow{3}{*}{\xmark} & \multirow{3}{*}{\checkmarkJTN}  \\
    of emergent ABM behaviour & & & \\
    throughout parameter space & & & \\
    \hline
    Extrapolates to unobserved  & \multirow{2}{*}{\checkmarkJTN} & \multirow{2}{*}{\checkmarkJTN} & \multirow{2}{*}{\xmark}  \\
    parameter values and data & & & \\
    \hline
    Interpretable & \xmark & \checkmarkJTN & \checkmarkJTN \\
    \hline
    
    \end{tabular}
    \caption{Highlighting the strengths (\checkmarkJTN) and limitations (\xmark) of the three approaches in this article for understanding the resulting behaviour from ABMs.}
    \label{tab:pro_cons}
\end{table}

\clearpage
\appendix

\section{Coarse graining BDM rules into a DE model} \label{sec:app_deriving_DE}

In this section, we will derive coarse-grained DE models of the BDM process. We begin by defining $\mathbb{P}[0_\alpha(t)]$ and $\mathbb{P}[A_\alpha(t)]$ as the probabilities that the individual lattice site $\alpha$ is either vacant or occupied, respectively, at time $t$. We simplify notation by writing: $\mathbb{P}[A_\alpha(t)]=C_\alpha(t)$ and $\mathbb{P}[0_\alpha(t)]=1-C_\alpha(t)$.

Similarly, let $\mathbb{P}[A_\alpha(t),A_\beta(t)]$ denote the probability that both neighbouring sites $\alpha$ and $\beta$ are occupied at time $t$; we refer to this value as the \emph{neighbouring lattice site occupancy probability}. Along these lines, $\mathbb{P}[0_\alpha(t),A_\beta(t)]$ is the probability that $\alpha$ is vacant and $\beta$ is occupied at time $t$, etc. These joint probabilities are related to the individual occupancy probabilities through their marginal probabilities:
\begin{align}
    C_\alpha(t) &= \mathbb{P}[A_\alpha(t),A_\beta(t)] + \mathbb{P}[A_\alpha(t),0_\beta(t)]; \nonumber \\
    1 - C_\alpha(t) &= \mathbb{P}[0_\alpha(t),A_\beta(t)] + \mathbb{P}[0_\alpha(t),0_\beta(t)]; \nonumber \\
    C_\beta(t) &= \mathbb{P}[A_\alpha(t),A_\beta(t)] + \mathbb{P}[0_\alpha(t),A_\beta(t)]; \nonumber \\
    1 - C_\beta(t) &= \mathbb{P}[A_\alpha(t),0_\beta(t)] + \mathbb{P}[0_\alpha(t),0_\beta(t)].
    \label{eq:marginal}
\end{align}
neighbouring occupancy probabilities are also related to the individual occupancy probabilities using the joint occupancy correlation function (see \citep{baker_correcting_2010}):
\begin{equation}
    F(t;\alpha,\beta) = \dfrac{\mathbb{P}[A_\alpha(t),A_\beta(t)]}{C_\alpha(t)C_\beta(t)}. \label{eq:correlation}
\end{equation}
Note that if $F(t;\alpha,\beta)=1$, then 
$\mathbb{P}[A_\alpha(t),A_\beta(t)] = C_\alpha(t)C_\beta(t),$
indicating that the occupancy of the neighbouring sites $\alpha$ and $\beta$ are independent. We can combine Equations \eqref{eq:marginal} and \eqref{eq:correlation} to write each neighbouring occupancy probability in terms of the individual occupancy probabilities and the occupancy correlation function 
\begin{align}
    \mathbb{P}[A_\alpha(t),A_\beta(t)] &= C_\alpha(t)C_\beta(t)F(t;\alpha,\beta); \nonumber \\
    \mathbb{P}[A_\alpha(t),0_\beta(t)] &= C_\alpha(t)\Big(1 - C_\beta(t)F(t;\alpha,\beta)\Big); \nonumber \\
    \mathbb{P}[0_\alpha(t),A_\beta(t)] &= C_\beta(t)\Big(1 - C_\alpha(t)F(t;\alpha,\beta)\Big); \nonumber \\
    \mathbb{P}[0_\alpha(t),0_\beta(t)] &= 1 - C_\alpha(t) - C_\beta(t) +   C_\alpha(t)C_\beta(t)F(t;\alpha,\beta). \label{eq:neighbors_interms_ind}
\end{align}

We are now ready to convert the rules of the BDM process into a coarse-grained DE model. We begin by writing a master equation for how $C_\alpha(t)$ will change due to the effects of agent birth, death, and migration:
\begin{equation}
    \dfrac{\text{d}C_\alpha(t)}{\text{d} t} = K_\text{birth} + K_\text{death} + K_\text{migration}.\label{eq:master_eqn_app}
\end{equation}
We now aim to derive espressions for $K_\text{birth}, K_\text{death},\text{ and } K_\text{migration}$. The birth reaction from Equation \eqref{eq:prolife_rate} specifies that the density at lattice site $\alpha$ may increase when $\alpha$ is unoccupied and $\beta\in\mathcal{B}(\alpha)$ is occupied because the agent at $\beta$ may give birth and place its daughter agent in $\alpha$. We can then write
\begin{equation}
    K_\text{birth} = \dfrac{P_p}{4} \sum_{\beta\in\mathcal{B}(\alpha)}  \mathbb{P}[0_\alpha(t),A_\beta(t)],
\end{equation}
because any of the neighbouring lattice sites may undergo birth events. Similarly, we can convert Equation \eqref{eq:death_rate} as
\begin{equation}
    K_\text{death} = -P_d C_\alpha(t),
\end{equation}
for agent death and convert Equation \eqref{eq:migrate_rate} as
\begin{equation}
    K_\text{migration} = \dfrac{P_m}{4}\sum_{\beta\in\mathcal{B}(\alpha)} \bigg( \mathbb{P}[0_\alpha(t),A_\beta(t)] - \mathbb{P}[A_\alpha(t),0_\beta(t)]\bigg),
\end{equation}
for agent migration. Substitution of these terms into Equation \eqref{eq:master_eqn_app} provides the master equation for the BDM process
\begin{align}
    \dfrac{\text{d}C_\alpha(t)}{\text{d} t} =  &\dfrac{P_p}{4} \sum_{\beta\in\mathcal{B}(\alpha)} \bigg( \mathbb{P}[0_\alpha(t),A_\beta(t)]\bigg) - P_d C_\alpha(t) + \dfrac{P_m}{4}\sum_{\beta\in\mathcal{B}(\alpha)} \bigg( \mathbb{P}[0_\alpha(t),A_\beta(t)] - \mathbb{P}[A_\alpha(t),0_\beta(t)]\bigg).  \label{eq:ABM_DE_nosimp_neighbor_pre}
\end{align}

Equation \eqref{eq:ABM_DE_nosimp_neighbor_pre} provides a DE model to describe the dynamics of $C_\alpha(t)$, however, this equation is not closed because we need to know $\mathbb{P}[A_\alpha(t),A_\beta(t)]$ in order to evaluate the right-hand side. We can use the marginal identities from Equations \eqref{eq:marginal} and \eqref{eq:neighbors_interms_ind} to simplify the terms in Equation \eqref{eq:ABM_DE_nosimp_neighbor_pre} and write
\begin{equation}
    \dfrac{\text{d}C_\alpha(t)}{\text{d} t} =   \dfrac{P_m}{4}\sum_{\beta\in\mathcal{B}(\alpha)} \bigg( C_\beta(t)  - C_\alpha(t) \bigg) + \dfrac{P_p}{4} \sum_{\beta\in\mathcal{B}(\alpha)} C_\beta(t) \bigg( 1 - C_\alpha(t) F(t;\alpha,\beta) \bigg) - P_d C_\alpha(t). \label{eq:ABM_DE_nosimp}
\end{equation}
We proceed by making a simplification in order to  close this system. Since we initiate all simulations with agents distributed uniformly at random, we assume that all individual occupancy probabilities are equally distributed\footnote{This assumption is not applicable when the initial configuration is spatially heterogeneous. See \citep{deutsch_cellular_2005,simpson_multi-species_2009} for the derivation of PDE models in this scenario.} on average so that $C_\alpha(t)=C_\gamma(t)=C(t)$ for any two lattice sites $\alpha$ and $\gamma$. From this assumption, we have that $F(t;\alpha,\beta)=F(t;|\alpha-\beta|)$, i.e., $F(t;\alpha,\beta)=F(t;1)$ for $\beta\in\mathcal{B}(\alpha)$. From Equation \eqref{eq:correlation}, we next write $\mathbb{P}[A_\alpha(t),A_\beta(t)]=C^2(t)F(t;1)$ for $\beta\in\mathcal{B}(\alpha)$. These observations lead to the following DE model from Equation \eqref{eq:ABM_DE_nosimp}:
\begin{equation}
    \dfrac{\text{d}}{\text{d}t}C(t) = P_pC(t)\big(1-C(t)F(t;1)\big) - P_dC(t). \label{eq:modified_logistic}
\end{equation}
Equation \eqref{eq:modified_logistic} is not yet closed because we still do not know $F(t;1)$. Our second simplification is the \emph{mean-field assumption}, in that the occupancy probabilities of neighbouring lattice sites are independent so that $F(t;1) \equiv 1$. This assumption leads to the \emph{mean-field model} for the ABM:
\begin{equation}
    \dfrac{\text{d}}{\text{d}t}C(t) = P_pC(t)\big(1-C(t)\big) - P_dC(t). \label{eq:MF_logistic_first_mention}
\end{equation}
Note that Equation \eqref{eq:MF_logistic_first_mention} can be re-formulated as the standard logistic DE model given by 
\begin{equation}
    \dfrac{\text{d}}{\text{d}t}C(t) = rC(t)\bigg(1-\dfrac{C(t)}{K}\bigg),
\end{equation}
where $r = P_p - P_d$,  $K = (P_p-P_d)/P_p$. This model is advantageous in that it is closed and can be solved analytically:
\begin{equation}
    C(t) = \dfrac{KC(0)\text{e}^{rt}}{K + C(0)(\text{e}^{rt}-1)},\label{eq:log_analytic}
\end{equation}
where $C(0)$ denotes the initial condition.

\clearpage
\section{Coarse graining SIR rules into a DE model} \label{sec:app_deriving_DE_SIR}

We now derive DE models governing the dynamics for $\mathbb{P}[S_\alpha(t)], \mathbb{P}[I_\alpha(t)],$ and $\mathbb{P}[R_\alpha(t)]$. As for the BDM model, because the initial agent configurations are uniformly distributed in space, we assume the probability of any type of agent occupancy ($S,I,$ or $R$) is independent of the lattice site and define $S(t) = \mathbb{P}[S_\alpha(t)], I(t) = \mathbb{P}[I_\alpha(t)], R(t) = \mathbb{P}[R_\alpha(t)].$ By converting the bimolecular reactions in Equations \eqref{eq:move_rate_SIR} and \eqref{eq:infect_rate} into the corresponding occupancy probability configurations that will lead to changes in $S,I,\text{ or }R$, and converting the monomolecular reaction in Equation \eqref{eq:recover_rate} into the individual occupancy probabilities that will lead to changes in $S,I,\text{ or }R$, we derive the master system of equations for $S(t),I(t),$ and $R(t)$ to be:
\begin{align}
    \dfrac{\text{d}}{\text{d}t}S(t) &= \sum_{\beta\in\mathcal{B}(\alpha)}\left[\dfrac{P_m}{4}\mathbb{P}[0_\alpha(t),S_\beta(t)] - \dfrac{P_m}{4}\mathbb{P}[S_\alpha(t),0_\beta(t)] - \dfrac{P_I}{4}\mathbb{P}[S_\alpha(t),I_\beta(t)] \right]; \nonumber \\
    \dfrac{\text{d}}{\text{d}t}I(t) &= \sum_{\beta\in\mathcal{B}(\alpha)}\left[\dfrac{P_m}{4}\mathbb{P}[0_\alpha(t),I_\beta(t)] - \dfrac{P_m}{4}\mathbb{P}[I_\alpha(t),0_\beta(t)] + \dfrac{P_I}{4}\mathbb{P}[I_\alpha(t),S_\beta(t)] \right] - R C_I(t); \nonumber \\
    \dfrac{\text{d}}{\text{d}t}R(t) &= \sum_{\beta\in\mathcal{B}(\alpha)}\left[\dfrac{P_m}{4}\mathbb{P}[0_\alpha(t),I_\beta(t)] - \dfrac{P_m}{4}\mathbb{P}[I_\alpha(t),0_\beta(t)]\right] + P_R I(t).\label{eq:SIR_DE_nosimp}
\end{align}
We then use the mean-field assumption to write $\mathbb{P}[Y_\alpha(t),Z_\beta(t)] = Y(t)Z(t)$, where $Y,Z\in\{S,I,R,0\}$. This assumption reduces Equation \eqref{eq:SIR_DE_nosimp} to the commonly used SIR model given by:
\begin{equation}
    \dfrac{\text{d}S}{\text{d}t} = -P_I SI; \ \ \ \ \  \dfrac{\text{d}I}{\text{d}t} = P_I SI - P_RI; \ \ \ \ \  \dfrac{\text{d}R}{\text{d}t} = P_RI. \label{eq:SIR_MF_app_prev}
\end{equation}
In Equation \eqref{eq:SIR_MF_app_prev}, the variables $S,I,$ and $R$ denote the density of susceptible, infected, and recovered agents over time, respectively, which cannot exceed 0.5 if only half of the simulation domain is occupied by agents. We can convert these variables to the fraction of susceptible, infected, and recovered agents by computing the dimensionless variables $S^*(t) = S(t)/M, I^*(t) = I(t)/M, \text{ and } R^*(t) = R(t)/M$, where $M$ is the proportion of occupied lattice sites in the simulation domain. The system of equations for these variables are given by
\begin{equation}
    \dfrac{\text{d}S^*}{\text{d}t} = -M P_I S^*I^*, \ \ \ \ \  \dfrac{\text{d}I^*}{\text{d}t} = M P_I S^*I^* - P_RI^*, \ \ \ \ \  \dfrac{\text{d}R^*}{\text{d}t} = P_RI^*. \label{eq:SIR_MF_app}
\end{equation}

\clearpage

\section{Gillespie algorithm}
\label{sec:Gillespie}
\begin{algorithm}
\SetAlgoLined
Create $X\times X$ lattice with user-specified placement of agents;

Maximum lattice occupancy is given by $N=X^2$

Set $t=0$; Set maximum simulation time $t_\text{end}$;

Set $C(t)$ equal to to number of agents on the lattice;

\While{$t<t_\text{end}$ and $C(t)< N$}{Randomly choose an agent and determine its lattice site;

Calculate the propensity function $a(t)=(P_p + P_m + P_d)C(t)$;

Calculate the following random variables, uniformly distributed on $[0,1]: \gamma_1,\gamma_2$;

Calculate time step $\tau=-\ln(\gamma_1)/a(t)$;

$t=t+\tau$;

$C(t)=C(t-\tau)$;

$R=a(t)\gamma_2$;

\uIf{$R<P_pC(t)$}{

Choose adjacent lattice site with equal probability 1/4;

\If{chosen lattice site is empty}{Place new agent to chosen lattice site;

$C(t) = C(t) + 1$;

}
}
\uElseIf{$R\in ( (P_pC(t),(P_p+P_m)C(t) ]$}{

Choose adjacent lattice site with equal probability 1/4;

\If{chosen lattice site is empty}{Move agent to chosen lattice site;}

}

\ElseIf{$R\in ( (P_p+P_m)C(t),(P_p+P_m+P_d)C(t) ]$}{

Remove agent from lattice site;

$C(t) = C(t) - 1$;
}

}

\caption{Gillespie algorithm for the BDM process (modified from \citep{fadai_accurate_2019}}\label{algo:gillespie}
\end{algorithm}

\clearpage

\section{Lasso algorithm using FISTA}
\label{sec:fista}
\begin{algorithm}
\SetAlgoLined
\SetKwInOut{input}{Input}
\SetKwInOut{output}{Output}

\textbf{Input:} Matrix $A$, vector $b$, and regularization parameter $\lambda$

\textbf{Output:} $x = \dfrac{1}{2}\arg\min\|Ax-b\|_2^2 + \lambda\|x\|_1.$

Get $d=$ \# columns in A.\\
\For{$i \gets 0$ \KwTo $d$ }{
        Set $a_i=\left\|A[:,i]\right\|_2$.\\
        Set $A[:,i] = A[:,i]/a_i.$  \\
    }
Set $L=\|A^TA\|_2$ (the largest singular value of $A$), $w = \vec{0}$, $w_{\text{old}} = \vec{0}$.\\

\While{$iter < iter_{\max}$}{
    $z = \dfrac{iter}{iter+1}(w-w_{\text{old}})$\\
    $w_\text{old} = w$\\
    $z = z - A^T(A^Tz-b)/L$\\
    \For{$i \gets 0$ \KwTo $d$ }{
        $w[i] = \text{sgn}(z[i])\times(\max\left\{\left|z[i]-\lambda/L \right|,0\right\})$\\
    }
}
\For{$i \gets 0$ \KwTo $d$ }{
        Set $w[i] = w[i]\times a_i.$\\
    }

\caption{Lasso implementation using FISTA from \citep{beck_fast_2009}.}\label{algo:lasso}
\end{algorithm}

\clearpage

\section{Hyperparameter selection for Lasso}
\label{sec:hyperselectlasso}
This section provides brief notes about the practical selection of hyperparameters for the Lasso method. We used regularization parameter $\lambda=0.0004$ for the Lasso algorithm to learn an equation for  $\mean{C_\text{ABM}(t)}$ in the tutorial in Section \ref{subsubsec:EQL_pipeline}. The optimal value of this hyperparameter is typically not known \emph{a priori}. There are several ways to select such a hyperparameter, including a grid search \citep{lagergren_learning_2020}, cross validation \citep{mangan2017model}, or Bayesian optimization \citep{snoek_practical_2012}. We discuss the grid search option here due to its simplicity. In a grid search to determine an appropriate value for $\lambda$, we specify several plausible options, given by $\{\lambda_1,\lambda_2,\dots,\lambda_n\}$ and split the data into  training (d$C_d^\text{train}(t)$/d$t$ and $\Theta^\text{train}$) and testing (d$C_d^\text{test}(t)$/d$t$ and $\Theta^\text{test}$) sets. The training and testing portions of $\Theta$ will contain all columns of $\Theta$ but only a subset of the rows. For a possible hyperparameter value, $\lambda_i$, we solve the Lasso problem from Equation \eqref{eq:lasso} using the training data and $\lambda=\lambda_i$ to determine the resulting $\xi$ estimates, $\hat{\xi}_i$. The optimal value of $\lambda$ is then chosen as:
\begin{equation}
    \hat{\lambda} = \arg\min_{\lambda_i} \left\| \dfrac{\text{d}C_d^\text{test}(t)}{\text{d}t} - \Theta^\text{test}\xi_i   \right\|_2.\label{eq:hyper_select}
\end{equation}

The optimal value $\hat{\lambda}$ results in the estimate $\hat{\xi}$ that best generalizes to the testing data. Recent work has shown that for Lasso, $\hat{\xi}$ tends to incorporate small additional terms in the final learned equation \citep{lagergren_learning_2020}. To select the regularization parameter, $\lambda$, for the Lasso algorithm, we randomly split half of $\mean{C_\text{ABM}(t)}$ into a training set and the remaining into a testing set. We considered 100 values of $\lambda$ between $10^{-5}$ and $10^{-3}$ as well as $\lambda=0$ and solved Equation \eqref{eq:lasso} using the training set for all 101 potential values of $\lambda$. We can perform this operation several times (changing the training and testing set each time) and notice that sometimes the chosen hyperparameter is zero and sometimes it is nonzero. When the chosen hyperparameter is zero, then the EQL pipeline learns an equation of the form
\begin{equation}
    \dfrac{\text{d}C}{\text{d}t} = 0.00488C    - 0.01187C^2    + 0.00806C^3    - 0.00831C^4,
\end{equation}
whereas when the chosen hyperparameter is nonzero, then the EQL pipeline learns an equation of the form
\begin{equation}
    \dfrac{\text{d}C}{\text{d}t} = 0.00468C    - 0.00951C^2.
\end{equation}

To ensure the final learned equation is sensitive to changes in each nonzero library coefficient, we perform a round of pruning after learning which proceeds as follows. The $j^{\text{th}}$ nonzero term of $\hat{\xi}$ is included in the final inferred equation if $\| \text{d}C_d^\text{test}(t)/\text{d}t - \Theta^\text{test}\hat{\xi}_j \|_2^2$ increases by a given pruning percentage, where here $\hat{\xi_j}$ is the estimated parameter vector with the $j^\text{th}$ term manually set to zero.  To ensure that the final learned equation is not an artifact of the training and testing split, we perform this entire process for ten randomized training and testing splits of the data and select the equation form arises most frequently. We set the parameters for each coefficient to be the mean of each coefficient for each time the final equation form was learned. If we set our pruning percentage to be 5\%, then the majority of learned equations will be of the form
\begin{equation}
    \dfrac{\text{d}C}{\text{d}t} = 0.00468C    - 0.00951C^2.
\end{equation}

\section{Case study 2: parameter distributions}\label{sec:bw}
\begin{figure}[h]
    \centering
    \includegraphics[width=0.89\textwidth]{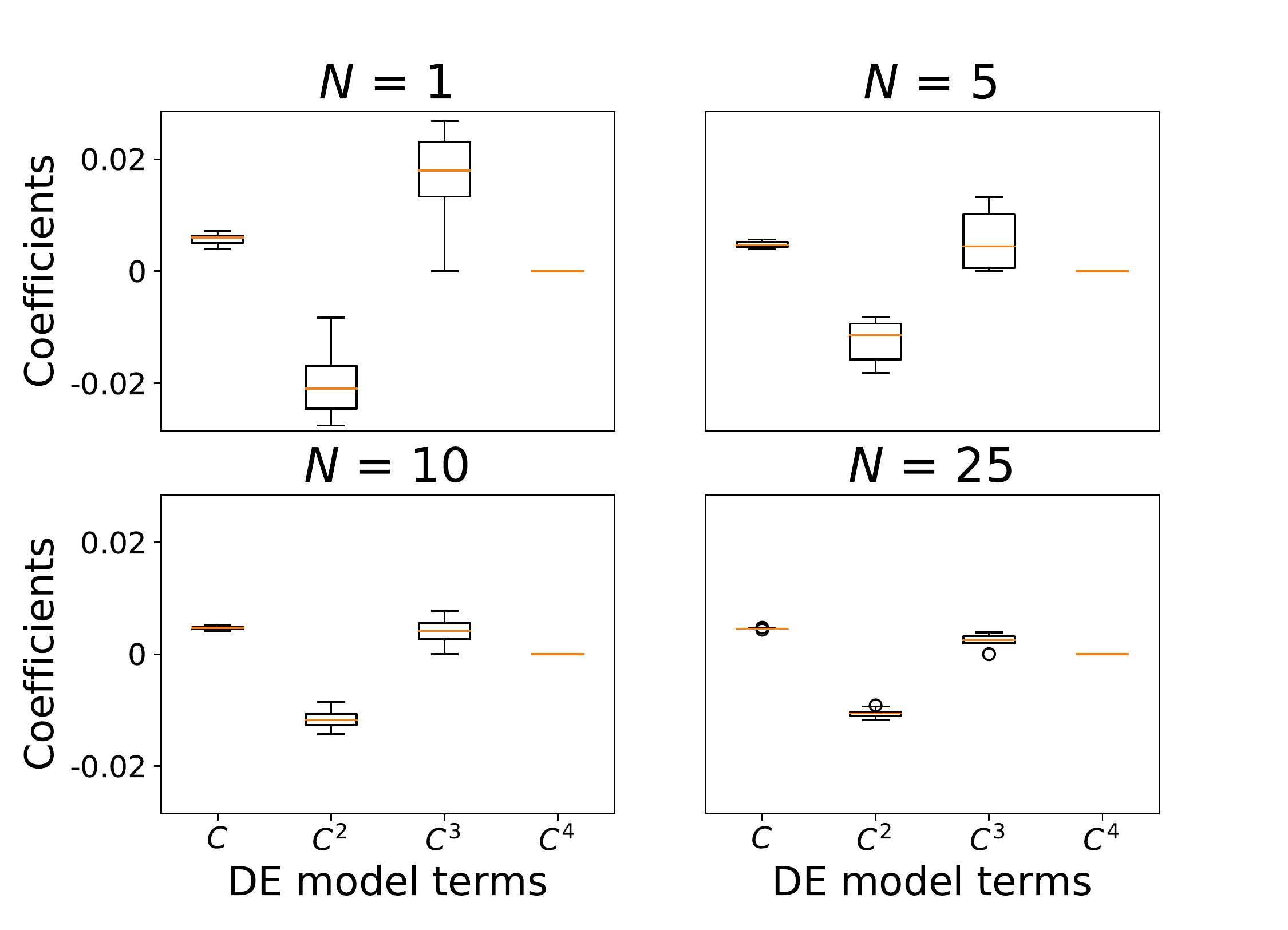}
    \caption{Summaries of the distributions of model parameters for the learned DE models from Case study 2 over various values of $N$. Each box and whisker plot summarizes the distribution of coefficient estimates from ten realisations of $\mean{C_\text{ABM}(t)}$ for various values of $N$. In each box and whisker plot, the lower line of the box portion provides the 25\% quartile of the data and the upper line denotes the 75\% quartile. The orange line on each box plot denotes the median coefficient value. The length of the upper and lower whiskers are 1.5 times the interquartile range of the distribution, and dots denote outlier points.}
    \label{fig:app_bw}
\end{figure}

\clearpage

\bibliographystyle{plainnat}      
\bibliography{references}


\end{document}